\documentclass[11pt,a4paper]{article}
\usepackage[a4paper,left=24mm,right=24mm,top=20mm,bottom=20mm]{geometry}
\usepackage{amsmath,amssymb,amsthm}
\usepackage{lmodern}
\usepackage[hidelinks]{hyperref}
\numberwithin{equation}{section}
\newcommand{\R}{\mathbb R}
\newcommand{\dd}{\,\mathrm d}
\newcommand{\HH}{\mathcal H}
\newcommand{\one}{\mathbf 1}
\newcommand{\PV}{\operatorname{PV}}
\newcommand{\dist}{\operatorname{dist}}
\newcommand{\Div}{\operatorname{div}}
\newcommand{\osc}{\operatorname{osc}}
\newtheorem{theorem}{Theorem}[section]
\newtheorem{lemma}[theorem]{Lemma}
\theoremstyle{definition}
\newtheorem{definition}{Definition}[section]
\hypersetup{pdftitle={Uniform interior gradient estimates for nonlocal minimal graphs},pdfauthor={Shuning Xu}}
\title{Uniform interior gradient estimates for nonlocal minimal graphs}
\author{Shuning Xu}
\date{}
\begin{document}
\maketitle
\begin{abstract}
For nonlocal minimal graphs of dimension $m\ge2$ we prove an interior gradient estimate whose exponent is linear in the height oscillation divided by the radius. Its constant stays bounded as $s\uparrow1$. A radial test of the fractional Jacobi inequality bounds a power of the reciprocal vertical normal at one point by a surface integral. Two tests of the graph equation and a squared-cutoff test of the Jacobi inequality bound that integral by the height oscillation.
\end{abstract}

\section{Introduction and main result}

A classical minimal graph $x_{m+1}=u(x')$ satisfies
\[
 \Div\!\left(\frac{Du}{\sqrt{1+|Du|^2}}\right)=0.
\]
An interior gradient estimate controls the slope of the graph away from the boundary in terms of its height. Finn obtained such estimates in dimension two \cite{Finn1953,Finn1954}. Finn \cite{Finn1963} and Serrin \cite{Serrin1963,Serrin1968} subsequently used comparison arguments to refine the bounds for positive solutions on a disk. In 1969, Bombieri, De Giorgi and Miranda proved an interior exponential gradient estimate in every dimension \cite{BDGM}. Their proof uses a Sobolev inequality on the graph and an iteration of De Giorgi--Moser type \cite{MirandaSobolev,SimonSurvey}. If $B'_r(x'_0)$ denotes the ball of radius $r$ centered at $x'_0$ in $\R^m$, the pointwise estimate takes the form
\[
 |Du(x'_0)|\le C_1(m)\exp\!\left[
 C_2(m)\frac{\sup_{B'_r(x'_0)}u-u(x'_0)}{r}
 \right].
\]

In 1972, Bombieri and Giusti proved a weak Harnack inequality on area-minimizing hypersurfaces by Moser iteration. Applying it to the vertical component of the unit normal, they obtained a polynomial gradient bound for entire minimal graphs \cite[Theorems~6 and~9 and (6.8)]{BG1972}. In the same year, Trudinger gave a direct proof of the local exponential estimate \cite[(19) and (22)--(28)]{Trudinger}. A radial potential reduces the pointwise estimate to an integral over the graph, and suitable tests of the minimal graph equation and the Jacobi inequality bound that integral in terms of the height. In a note added in proof, Trudinger credits Michael with a different proof of the mean-value inequality and Simon with its extension to arbitrary codimension \cite{Trudinger,MichaelSimon1973}. Korevaar and Wang later gave maximum-principle arguments applicable to the minimal graph equation \cite{Korevaar1986,Wang1998}. De Silva and Jerison obtained another proof using a normal deformation of the graph and a weak Harnack inequality \cite[Theorem~1.2 and Section~3]{DeSilvaJerison}.

We follow Cabr\'e and Cozzi \cite[Section~1.1]{CC} in defining the fractional perimeter. Fix $0<s<1$. For measurable sets $A,B,E\subset\R^{m+1}$ with $A\cap B=\varnothing$ and a bounded open set $V\subset\R^{m+1}$, we set
\[
\begin{gathered}
 I_s(A,B)=\int_A\int_B\frac{\dd y\,\dd x}{|x-y|^{m+1+s}},\\
 \operatorname{Per}_s(E;V)=I_s(E\cap V,E^c)+I_s(E\setminus V,V\setminus E).
\end{gathered}
\]
We say that $E$ is $s$-minimal in $V$ if $\operatorname{Per}_s(E;V)$ is finite and $\operatorname{Per}_s(E;V)\le\operatorname{Per}_s(F;V)$ for every measurable $F$ such that $F\setminus V=E\setminus V$ almost everywhere.

For an open set $U$, we call $E$ locally $s$-minimal in $U$ if it is $s$-minimal in every open set $V\Subset U$. We then call $\partial E$ an $s$-minimal surface in $U$ \cite{CRS,CC,CL}.

At a point $x\in\partial E$ where the boundary is $C^2$ in a neighborhood, the fractional mean curvature of $E$ is defined by
\[
 H_s[E](x)
 =\PV\int_{\R^{m+1}}
       \frac{\one_{E^c}(y)-\one_E(y)}{|x-y|^{m+1+s}}\,\dd y .
\]
The principal value is taken by excluding $B_\varepsilon(x)$ and letting $\varepsilon\downarrow0$. If $E$ is locally $s$-minimal in $U$, the first variation of the fractional perimeter gives the Euler--Lagrange equation $H_s[E](x)=0$ at every such point in $U$ \cite[(1.4)--(1.5)]{CC}.

Caffarelli, Roquejoffre and Savin established the variational theory and the first interior regularity results for these surfaces in their 2010 paper \cite{CRS}. Ambrosio, De Philippis and Martinazzi proved the $\Gamma$-convergence of the normalized fractional perimeter and a convergence theorem for local minimizers in 2011 \cite{ADPM}. Caffarelli and Valdinoci obtained density estimates uniform as $s\uparrow1$ in 2011 and an improvement of flatness uniform in the same limit in 2013 \cite{CV11,CV13}. These results make it possible to use the classical theory of perimeter minimizers to study fractional minimizers for $s$ close to one.

Barrios, Figalli and Valdinoci proved in 2014 that $s$-minimal boundaries with suitable $C^{1,\alpha}$ regularity are smooth \cite[Theorem~1.1]{BFV2014}. In 2016, Dipierro, Savin and Valdinoci proved that, under suitable assumptions on the cylinder and the exterior graph data, the minimizing set is a subgraph \cite[Theorem~1.1]{DSV2016}. Figalli and Valdinoci showed that locally Lipschitz $s$-minimal graphs are smooth in the interior \cite[Theorem~1.1]{FV}. Their paper appeared in print in 2017. This left the quantitative question of how to bound the gradient of a nonlocal minimal graph using only its local height oscillation.

Cabr\'e and Cozzi answered this question in a paper first circulated in 2017 and published in 2019. Suppose that the global subgraph of $u$ is locally $s$-minimal in $B'_{2r}\times\R$ and that $u$ is bounded in $B'_r$, where $B'_r=B'_r(0)$. In our notation their estimate is
\[
 \sup_{B'_r}|Du|
 \le C(m,s)\left(1+\frac{\osc_{B'_r}u}{r}\right)^{m+1+s}.
\]
They apply a weak Harnack inequality to the vertical component of the unit normal, which satisfies a truncated fractional Jacobi inequality \cite[Theorems~1.1 and~1.3(ii)]{CC}. The proof of the weak Harnack inequality combines fractional Sobolev inequalities on the surface with Moser iteration \cite[Theorem~1.5, Corollary~5.3 and Theorem~6.4]{CC}. Fall subsequently combined their gradient bound with higher regularity estimates to control interior derivatives of all orders in terms of the height oscillation \cite[Theorem~1.5 and the discussion following it]{Fall}.

Cabr\'e and Cozzi also observed that for fixed $s$ one can obtain exponential dependence on the oscillation by a covering argument based on their weak Harnack inequality. They pointed out that the sharp exponential growth in the classical problem forces the constant in their polynomial estimate to diverge as $s\uparrow1$. They raised the question of obtaining a gradient bound uniform in this limit and suggested starting from a fractional Sobolev inequality on the graph with the appropriate dependence on $s$ \cite[discussion preceding and following (1.2)]{CC}. We adapt Trudinger's direct argument using a radial potential and test functions. The proof yields a uniform exponential gradient bound without a fractional Sobolev inequality on the graph.

\begin{theorem}\label{thm:gradient}
For each $m\ge2$ there exist $s_*(m)\in(0,1)$ and $C(m)>0$ depending only on $m$. Let $s_*(m)\le s<1$ and $R>0$. Let $E=\{(x',t)\in\R^m\times\R:t<u(x')\}$ be the global subgraph of a measurable function $u:\R^m\to\R$. If $E$ is locally $s$-minimal in $B'_{4R}\times\R$ and $u$ is smooth with finite oscillation in $B'_{2R}$, then
\begin{equation}\label{eq:main}
 \sup_{B'_R}\sqrt{1+|Du|^2}
 \le \exp\!\left[C(m)\left(1+\frac{\osc_{B'_{2R}}u}{R}\right)\right].
\end{equation}
\end{theorem}

For an essentially bounded graph function one can take the smooth representative given by Lemma~\ref{lem:smooth-representative}. If $s_j\uparrow1$ and the oscillations of the corresponding functions $u_j$ in $B'_{2R}$ are uniformly bounded, we can vertically translate the graphs so that $u_j(0)=0$. Estimate \eqref{eq:main} gives local uniform convergence of a subsequence in $B'_R$ to a Lipschitz function. Local uniform convergence identifies the subgraph of this function with the local $L^1$ limit in $B'_R\times\R$, so it minimizes classical perimeter there by Lemma~\ref{lem:geometry}.

Section~\ref{sec:prelim} establishes local area estimates uniform as $s\uparrow1$ and the truncated Jacobi inequality used in the proof. In Section~\ref{sec:proof}, a radial potential reduces the pointwise bound to a surface integral, which we estimate in terms of the height oscillation by testing the equations.

\paragraph{Acknowledgments.} The author thanks Professor Xi-Nan Ma for proposing this problem. This work was supported by the National Key R\&D Program of China (2025YFA1017603).

\section{Preliminaries}\label{sec:prelim}
\subsection{Definitions and preliminary results}
For balls and cylinders we follow the conventions of Cabr\'e and Cozzi \cite[Section~2]{CC}. We write $m$ for the graph dimension and $s$ for the perimeter parameter.

\begin{definition}
For $P\in\R^{m+1}$, write $P'\in\R^m$ for its horizontal projection. For $r>0$ the ambient ball, base ball, and vertical cylinder are
\[
\begin{aligned}
 B_r(P)&=\{y\in\R^{m+1}:|y-P|<r\},\\
 B'_r(P')&=\{y'\in\R^m:|y'-P'|<r\},\\
 C_r(P')&=B'_r(P')\times\R.
\end{aligned}
\]
When the center is omitted, balls are centered at the origin and cylinders are centered over the origin.
\end{definition}

\begin{definition}
For $E\subset\R^{m+1}$ measurable and $U\subset\R^{m+1}$ open, the relative classical perimeter is
\[
\operatorname{Per}(E;U)
 =\sup_{\substack{\varphi\in C_c^1(U;\R^{m+1})\\|\varphi|\le1}}
     \int_U\one_E\,\Div\varphi\,\dd x.
\]
When $E$ has locally finite perimeter, this is $|D\one_E|(U)$, the total variation in $U$ of the distributional derivative of its indicator.
\end{definition}

We choose a representative of $E$ that contains every point with a neighborhood in which $E^c$ has measure zero and excludes every point with a neighborhood in which $E$ has measure zero. Thus
\[
 P\in\partial E
 \quad\Longleftrightarrow\quad
 0<|E\cap B_r(P)|<|B_r(P)|\quad\text{for every }r>0.
\]
Write $\Sigma=\partial E$. On its smooth part, $\nu$ is the exterior unit normal and $\dd\sigma=\dd\HH^m$ is surface measure, as in Cabr\'e and Cozzi. For a set of locally finite perimeter the area measure is $\HH^m\llcorner\partial^*E$, where $\partial^*E$ is the reduced boundary. On the graph of $u$ the vertical component is $\omega=\nu_{m+1}=(1+|Du|^2)^{-1/2}$. In expressions on the base, $f(x')$ denotes the value of the surface function $f$ at $(x',u(x'))$. The graph area formula is $\omega\dd\sigma=\dd x'$.

\begin{definition}
Write $K(x,y)=|x-y|^{-m-1-s}$. For a smooth interior point $x$ of an integration patch $D\subset\Sigma$ and a function $f$ smooth near $x$ with integrable tail, set
\begin{equation}\label{eq:operators}
 L_D f(x)=(1-s)\PV\int_D(f(x)-f(y))K(x,y)\dd\sigma(y).
\end{equation}
The principal value is taken by deleting a ball centered at $x$ and shrinking its radius to zero. The operator has order $1+s$. For $D=\Sigma\cap U$, it is the operator in \cite[(1.15)]{CC} multiplied by $-(1-s)$.
\end{definition}

Throughout the estimates we take $s\ge3/4$ and raise its lower threshold when needed. Unless indicated otherwise, $c,C>0$ depend only on the graph dimension $m$ and may change from line to line. The lower threshold for $s$ also depends only on $m$. Other dependencies in arguments at fixed $s$ or on a fixed surface are indicated where used.

Cinti, Serra and Valdinoci proved a local $BV$ bound for stable sets \cite[Theorem~1.1]{CSV}. An $s$-minimal set is stable, so their bound gives the following estimate after rescaling.

\begin{lemma}\label{lem:fixed-perimeter}
Fix $s\in(0,1)$ and $r>0$. If $E$ is $s$-minimal in $B_{2r}(P)$, then its classical perimeter satisfies
\[
 \operatorname{Per}(E;B_r(P))\le C(m,s)r^m.
\]
The constant is independent of $E,P,r$.
\end{lemma}

At $r=1$ it makes the supremum of $\operatorname{Per}(E;B_1)$ over $s$-minimal sets in $B_4$ finite at each fixed $s$. Lemma~\ref{lem:geometry} obtains a bound for this supremum as $s\uparrow1$.

Caffarelli and Valdinoci established density estimates uniform as $s\uparrow1$ \cite{CV11}. We apply the estimate to $E$ and to its complement to obtain both bounds.

\begin{lemma}\label{lem:density}
There are $s_*(m)\in(0,1)$ and $c(m)>0$ depending only on $m$. If $s\in[s_*(m),1)$, $P\in\partial E$, and $E$ is $s$-minimal in $B_{2r}(P)$ for $r>0$, then
\[
|E\cap B_r(P)|\ge c(m)r^{m+1},
\qquad |B_r(P)\setminus E|\ge c(m)r^{m+1}.
\]
\end{lemma}

We use the following form of the improvement of flatness of Caffarelli and Valdinoci \cite[Theorem~1 and its proof in Sections~2.6--2.7]{CV13}.

\begin{lemma}\label{lem:flat}
Fix $s_0,\gamma\in(0,1)$. There are $\delta_0,\vartheta\in(0,1)$ depending only on $m,s_0,\gamma$. Let $s\in[s_0,1)$, $r>0$, $P\in\partial E$, and suppose $E$ is $s$-minimal in $B_r(P)$. Suppose some unit vector $e\in\mathbb S^m$ satisfies the slab condition
\begin{equation}\label{eq:flatness}
\partial E\cap B_r(P)
\subset\{x\in\R^{m+1}: |(x-P)\cdot e|\le\delta_0 r\}.
\end{equation}
Then in $B_{\vartheta r}(P)$ the boundary consists of a single graph over the plane perpendicular to $e$, with height measured along $e$. In the coordinates obtained by translating $P$ to the origin and rescaling lengths by $r$, the graph has a $C^{1,\gamma}$ bound depending only on $m,s_0,\gamma$.
\end{lemma}

The next lemma uses Theorems~1 and~2 of Ambrosio, De Philippis and Martinazzi for compactness and the lower bound, and the recovery comparison in their proof of Theorem~3 for minimality \cite{ADPM}. A truncation argument adapts the comparison to local convergence in $U$.

\begin{lemma}\label{lem:compactness}
Let $s_j\uparrow1$ and let $U\subset\R^{m+1}$ be open. Assume that for every ball $B\Subset U$,
\[
 \sup_j (1-s_j)\int_{E_j\cap B}\!\int_{B\setminus E_j}
           \frac{\dd y\dd x}{|x-y|^{m+1+s_j}}<\infty.
\]
Both endpoints of each pair in this integral lie in $B$.

The indicators $\one_{E_j}$ are relatively compact in $L^1_{\rm loc}(U)$. If $\one_{E_j}\to\one_{E_\infty}$ along a subsequence, then $E_\infty$ has locally finite perimeter in $U$, and for every $B\Subset U$ along that subsequence,
\[
 |B'_1|\operatorname{Per}(E_\infty;B)
 \le\liminf_{j\to\infty}(1-s_j)
      \int_{E_j\cap B}\!\int_{B\setminus E_j}
      \frac{\dd y\dd x}{|x-y|^{m+1+s_j}}.
\]
If each $E_j$ is locally $s_j$-minimal in $U$, then $E_\infty$ locally minimizes classical perimeter in $U$.
\end{lemma}

\begin{proof}
The compactness and the lower bound follow by local application of the cited results. To prove minimality, take bounded smooth open sets $V\Subset W\Subset U$, truncate $E_j$ and $E_\infty$ to $W$, so their indicators converge in $L^1(\R^{m+1})$. For competitors agreeing with $E_j$ outside $V$, this truncation changes the difference of the normalized energies by at most
\[
 2(1-s_j)\int_V\int_{W^c}|x-y|^{-m-1-s_j}\dd y\dd x
 \le \frac{C(m)(1-s_j)|V|}{s_j\,\dist(V,W^c)^{s_j}}\longrightarrow0 .
\]
Using this bound in the recovery comparison and passing to the limit gives the classical perimeter comparison. Since $V\Subset U$ is arbitrary, the limit locally minimizes perimeter in $U$.
\end{proof}

The following property of the singular set is a consequence of the classical dimension estimate \cite[Chapter~7, Theorem~5.8]{Simon}.

\begin{lemma}\label{lem:singular}
Let $E$ locally minimize classical perimeter in an open set $U$. Let $S$ be the set of points of $\partial E\cap U$ that have no neighborhood in which the boundary is a smooth hypersurface. Then
\[
 \HH^{m-1}(S)=0.
\]
\end{lemma}

For fixed $s$, the existence of a smooth representative follows by applying Cabr\'e and Cozzi \cite[Theorem~1.1]{CC} on interior balls. It also follows from the implication (iv) $\Rightarrow$ (v) of \cite[Theorem~1.10]{CL}.

\begin{lemma}\label{lem:smooth-representative}
Let $u:\R^m\to\R$ be measurable with $u\in L^\infty(B'_{2R})$ for $R>0$. If its global subgraph is locally $s$-minimal in $B'_{4R}\times\R$, then $u$ agrees almost everywhere on $B'_{2R}$ with a smooth function. With the representative of $E$ chosen above, $\partial E\cap C_{2R}$ is the graph of this function.
\end{lemma}

\subsection{Uniform area and the truncated Jacobi inequality}
We use Lemma~\ref{lem:flat} from here on with $s_0=3/4$ and $\gamma=1/2$. The annular sums in Section~\ref{sec:proof} require $\gamma<s$ and $2\gamma<1+s$. This choice leaves a fixed positive margin in both inequalities for $s\ge3/4$.

We combine compactness, the density bounds and improvement of flatness to derive a classical perimeter bound uniform as $s\uparrow1$. The regular charts control the perimeter away from the singular set, which is handled using Lemma~\ref{lem:singular}.

\begin{lemma}\label{lem:geometry}
There is a dimensional threshold $s_*(m)\in[3/4,1)$ with the following properties. If $s\in[s_*(m),1)$ and $E$ is $s$-minimal in $B_{4r}(P)$, then
\begin{equation}\label{eq:area}
 \operatorname{Per}(E;B_r(P))\le C(m)r^m.
\end{equation}
If $P\in\partial E$, then $\operatorname{Per}(E;B_r(P))\ge c(m)r^m$.

If $s_j\uparrow1$ and each $E_j$ is locally $s_j$-minimal in a common open set $U$, a subsequence satisfies $\one_{E_j}\to\one_{E_\infty}$ in $L^1_{\rm loc}(U)$ for a classical perimeter minimizer $E_\infty$. Along this subsequence the boundaries converge locally in Hausdorff distance, with local $C^1$ convergence as single graphs near each regular point of $\partial E_\infty$. The classical perimeter measures converge weakly locally to $|D\one_{E_\infty}|$.
\end{lemma}

\begin{proof}
We first bound the normalized fractional perimeters in order to apply Lemma~\ref{lem:compactness}. Write $B_r=B_r(P)$. The competitor $E\setminus B_r$ gives
\[
\operatorname{Per}_s(E;B_r)
\le \int_{B_r}\int_{B_r^c} K(x,y)\dd y\dd x
\le \frac{C(m)r^{m+1-s}}{s(1-s)} .
\]
Indeed, the inner integral is bounded by $C(m)\dist(x,\partial B_r)^{-s}/s$. Integrating the distance to the sphere gives the stated bound. For $s_j\uparrow1$ and $E_j$ locally $s_j$-minimal in $U$, this yields
\[
(1-s_j)\int_{E_j\cap B_r}\int_{B_r\setminus E_j}
    |x-y|^{-m-1-s_j}\dd y\dd x
\le (1-s_j)\operatorname{Per}_{s_j}(E_j;B_r)
\le C(m)r^{m+1-s_j}
\]
whenever $B_r\Subset U$. The first integral counts only pairs with both endpoints in the ball. Lemma~\ref{lem:compactness} now gives $\one_{E_j}\to\one_{E_\infty}$ in $L^1_{\rm loc}(U)$ along a subsequence, classical minimality of $E_\infty$, and
\[
|B'_1|\operatorname{Per}(E_\infty;B_r)
\le \liminf_j (1-s_j)\operatorname{Per}_{s_j}(E_j;B_r)
\le C(m)r^m.
\]
The uniform classical upper bound for the approximating sets will follow below by excluding area concentration.

For the lower perimeter bound, we use the relative isoperimetric inequality on $B_r$ \cite[Section~5.6.2]{EG}, which gives
\[
\min\{|E\cap B_r|,|B_r\setminus E|\}^{m/(m+1)}
\le C(m)\operatorname{Per}(E;B_r).
\]
If $P\in\partial E$, Lemma~\ref{lem:density} bounds both volumes from below by $c(m)r^{m+1}$, so the left side is at least $c(m)r^m$.

We next prove local Hausdorff convergence of the boundaries. If a ball $B\Subset U$ is centered on $\partial E_\infty$, our choice of representative gives $0<|E_\infty\cap B|<|B|$. Hence
\[
\min\{|E_j\cap B|,|B\setminus E_j|\}
\ge\min\{|E_\infty\cap B|,|B\setminus E_\infty|\}
   -|(E_j\triangle E_\infty)\cap B|>0
\]
for large $j$. Thus $\partial E_j$ meets $B$. Conversely a limiting point of $\partial E_j$ outside $\partial E_\infty$ would lie in a ball in which $E_\infty$ or its complement is null. Lemma~\ref{lem:density} on a smaller ball about the approximating points contradicts $L^1$ convergence. By compactness, for each compact $F\Subset U$,
\[
\sup_{x\in F\cap\partial E_j}\dist(x,\partial E_\infty)
+\sup_{x\in F\cap\partial E_\infty}\dist(x,\partial E_j)\longrightarrow0
\]
with empty suprema taken as zero.

Around a regular point $P$ of $\partial E_\infty$, smoothness gives a radius $r$ such that the boundary in $B_{2r}(P)$ lies within $\delta_0 r/2$ of its tangent plane. Choose $P_j\in\partial E_j$ tending to $P$. Hausdorff convergence, with the spare half of the width $\delta_0r$, transfers \eqref{eq:flatness} to $B_r(P_j)$ with the same normal. Lemma~\ref{lem:flat} gives single graphs on a common smaller disk with uniform $C^{1,\gamma}$ bounds. Compactness in $C^1$ and Hausdorff convergence then identify the limit graph with the regular patch of $\partial E_\infty$.

We now prove the uniform classical perimeter bound \eqref{eq:area} by contradiction. At unit scale let $A_s$ be the supremum of $\operatorname{Per}(E;B_1)$ over $s$-minimal sets in $B_4$. It is finite at fixed $s$ by Lemma~\ref{lem:fixed-perimeter}. If it were unbounded as $s\uparrow1$, choose $s_j\uparrow1$ and $E_j$ with $A_{s_j}\to\infty$ and $\operatorname{Per}(E_j;B_1)\ge A_{s_j}/2$. The compactness just proved applies in $B_4$. For $x\in B_2$ and $0<r<1/4$, translate and dilate $B_{4r}(x)$ to $B_4$. The classical perimeter scales by $r^m$, so
\[
\frac{\operatorname{Per}(E_j;B_r(x))}{A_{s_j}}\le r^m.
\]
Write $\mu$ for a local weak limit of $|D\one_{E_j}|/A_{s_j}$ in $B_2$. A cutoff equal to one on $B_1$ and supported in $B_{3/2}$ gives $\mu(B_{3/2})\ge1/2$. Off $\partial E_\infty$ the measures vanish locally. At a regular point choose a fixed graph ball $B$; the $C^1$ convergence gives
\[
0\le\frac{\operatorname{Per}(E_j;B)}{A_{s_j}}\le\frac{C(B)}{A_{s_j}}\longrightarrow0.
\]
Thus $\mu$ is supported on $S$. The ball bound passes to the limit: $\mu(B_t(z))\le t^m$ for $z\in B_2$, $B_t(z)\Subset B_2$, $t<1/4$. Lemma~\ref{lem:singular} gives $\HH^m(S)=0$. Covering $S\cap B_{3/2}$ by such balls gives $\mu(B_{3/2})\le\sum_i t_i^m$ with arbitrarily small sum, a contradiction.

Apply the resulting local bound to the unnormalized perimeter measures. Every weak limit has zero mass on $S$ by the same growth argument. In the common regular charts the area formula and $C^1$ convergence identify the limit with $|D\one_{E_\infty}|$ \cite[Section~3.3]{EG}. This proves the asserted weak convergence.
\end{proof}

The compactness in Lemma~\ref{lem:geometry} also lets us derive a lower bound for the energy of the normal when the slab condition \eqref{eq:flatness} fails.

\begin{lemma}\label{lem:nonflat}
Take $\delta_0$ from Lemma~\ref{lem:flat} with the parameters fixed above. There exist $A=A(m)\ge10$, $c=c(m)>0$ and a dimensional threshold $s_*$. Suppose $P\in\partial E$, $E$ is $s$-minimal in $B_{4Ar}(P)$, and the boundary is smooth in $B_{Ar}(P)$. If for every $e\in\mathbb S^m$ there is $x\in\partial E\cap B_r(P)$ with $|(x-P)\cdot e|>\delta_0r$, then
\begin{equation}\label{eq:nonflat}
(1-s)\iint_{(\Sigma\cap B_{Ar}(P))^2}
\frac{1-\nu(x)\cdot\nu(y)}{|x-y|^{m+1+s}}\dd\sigma(x)\dd\sigma(y)
\ge c r^{m-1-s}\qquad(s\ge s_*).
\end{equation}
\end{lemma}

\begin{proof}
Translate and scale to $P=0$ and $r=1$. If the assertion fails, for $A_j=j+10$ we may choose $s_j\uparrow1$ and sets $E_j$ minimizing in $B_{4A_j}$, with $0\in\partial E_j$ and smooth boundary in $B_{A_j}$, such that the integral in \eqref{eq:nonflat}, including its factor $1-s_j$, is less than $1/j$. For each $j$, no unit vector satisfies \eqref{eq:flatness} in $B_1$. For any fixed ball $B_i$, $B_i\Subset B_{4A_j}$ for large $j$, so Lemma~\ref{lem:geometry} applies there. A diagonal subsequence gives a classical limit $E_\infty$ minimizing perimeter under compactly supported changes in every ball. Its boundary contains $0$ by Hausdorff convergence.

Near a regular limit point choose a common graph disk $B'$. Restrict both points in \eqref{eq:nonflat} to the corresponding graph patch, using the nonnegative integrand. In the graph coordinates the area factors and chord lengths are uniformly comparable to the planar ones. Writing $\nu_j$ for the pulled-back normals, the identity $2(1-\nu_j(x)\cdot\nu_j(y))=|\nu_j(x)-\nu_j(y)|^2$ gives

\[
(1-s_j)\iint_{B'\times B'}\frac{|\nu_j(x)-\nu_j(y)|^2}{|x-y|^{m+1+s_j}}\dd x\dd y\longrightarrow0 .
\]

For $\varphi\in C_c^\infty(B')$, interchanging $x,y$ in the truncated integrals gives, componentwise,

\[
\int_{B'}\nu_j L_{B'}\varphi\dd x
=\frac{1-s_j}{2}\iint_{B'\times B'}
  \frac{(\nu_j(x)-\nu_j(y))(\varphi(x)-\varphi(y))}
  {|x-y|^{m+1+s_j}}\dd x\dd y .
\]

The normalized squared-difference integral for $\varphi$ is uniformly bounded: its local radial integral is bounded by a constant times $(1-s_j)\int_0^{\operatorname{diam}B'} t^{-s_j}\dd t$. Cauchy--Schwarz and the vanishing normal energy make the pairing tend to zero. To compute $L_{B'}\varphi$, take $\ell>0$ less than one third of the distance from the support to the boundary. For points within $\ell$ of the support the entire ball of radius $\ell$ lies in $B'$. The linear Taylor term integrates to zero and the quadratic term gives

\[
-\frac{|B'_1|}{2}\Delta\varphi(x)
\,(1-s_j)\int_0^\ell t^{-s_j}\dd t
=-\frac{|B'_1|}{2}\ell^{1-s_j}\Delta\varphi(x),
\]

where the angular integral of the product of two coordinate components of the unit radial vector is $|B'_1|$ for equal indices and zero otherwise. The Taylor remainder is bounded by $C(1-s_j)\int_0^\ell t^{1-s_j}\dd t$. The tail is bounded by $C(1-s_j)\int_\ell^\infty t^{-2-s_j}\dd t$, which controls the entire integral at points separated from the support by at least $\ell$. Thus the operator converges uniformly to $-(|B'_1|/2)\Delta\varphi$. Passing to the limit in the pairing shows $\Delta\nu_\infty=0$ on the disk. Since $|\nu_\infty|=1$, $\Delta|\nu_\infty|^2=2|D\nu_\infty|^2$ makes the normal constant. The regular patch is therefore open in an affine plane.

Let $\Pi$ be the plane of such a patch and $S$ the singular set. Since $\HH^{m-1}(S)=0$, $\Pi\setminus S$ is path connected. For two of its points, project $S\cap \Pi$ onto the $(m-1)$-plane normal to the segment joining the points, within $\Pi$. The projection is null. A small parallel translate of the segment avoids $S$, and its endpoints can be joined to the original endpoints in neighborhoods disjoint from $S$. The points where the boundary locally equals $\Pi$ form a relatively open and closed subset of $\Pi\setminus S$. For closedness, a limit point is on the boundary and is regular, and its planar patch shares an open piece with $\Pi$. Taking the closure gives $\Pi\subset\partial E_\infty$. Distinct such planes are parallel, since an intersection would lie in $S$ and have dimension $m-1$. Each plane meeting a bounded ball contributes a fixed positive amount of perimeter in a larger concentric ball. Thus the planes are locally finite, and their union is the entire boundary by density of the regular points.

If two adjacent planes are at $t=0$ and $t=d$, replace the phase in $B'_L\times(0,d)$ by the phase just outside the two planes. In a ball containing the closed cylinder, the change in perimeter is

\[
m|B'_1|L^{m-1}d-2|B'_1|L^m.
\]

The first term is the new cylindrical side and the second accounts for the two disks removed. For $L>md/2$ the change is negative, contrary to minimality. Thus the boundary is a single plane through $0$. If $e$ is its unit normal, $|X\cdot e|=\dist(X,\partial E_\infty)$. By local Hausdorff convergence this distance is uniformly less than $\delta_0$ on $\partial E_j\cap B_1$ for large $j$. This yields \eqref{eq:flatness} in direction $e$, a contradiction.
\end{proof}

Translate the base center to the origin and assume that $C_{2r}$ lies in the smooth minimizing region. The following area bound holds for every ambient ball:

\begin{equation}\label{eq:all-growth}
\HH^m(\Sigma\cap C_r\cap B_t(x))\le C(m)t^m .
\end{equation}

For $t<r/10$ with nonempty intersection, take $P\in\Sigma\cap C_r\cap B_t(x)$. Then $B_t(x)\subset B_{2t}(P)$ and $B_{8t}(P)\subset C_{2r}$, so \eqref{eq:area} gives the bound. For $t\ge r/10$, let $\mathcal Z$ be the cylinder with base $B'_{r+r/100}$ and vertical interval $[x_{m+1}-t-r/100,x_{m+1}+t+r/100]$. For every $y\in\Sigma\cap C_r\cap B_t(x)$, the centers for which $y\in B_{r/100}(z)$ form exactly $B_{r/100}(y)$, contained in $\mathcal Z$. Fubini's theorem gives

\[
\int_{\mathcal Z}\HH^m(\Sigma\cap C_r\cap B_t(x)\cap B_{r/100}(z))\dd z
=|B_{r/100}|\,\HH^m(\Sigma\cap C_r\cap B_t(x)).
\]

For each such center, $B_{4r/100}(z)\subset C_{2r}$, so \eqref{eq:area} bounds the integrand by $Cr^m$. The volume of $\mathcal Z$ is at most $Cr^m(t+r)$. Dividing by $|B_{r/100}|$ gives $Cr^{m-1}(t+r)\le Ct^m$. To estimate the total area of $\Sigma\cap C_r$ when $\osc_{B'_r}u$ is finite, use the full height interval of the graph instead, enlarging it by $r/100$ at each end so that each required $B_{r/100}(y)$ remains inside the region of centers. This gives $\HH^m(\Sigma\cap C_r)\le Cr^{m-1}(r+\osc_{B'_r}u)$.

The growth bound \eqref{eq:all-growth} controls the mass of every radial layer. For all $\ell>0$ it gives

\begin{equation}\label{eq:moments-upper}
\begin{aligned}
(1-s)\int_{\Sigma\cap C_r\cap B_\ell(x)} |x-y|^2 K(x,y)\dd\sigma(y)
 &\le C\ell^{1-s},\\
\int_{(\Sigma\cap C_r)\setminus B_\ell(x)} |x-y|K(x,y)\dd\sigma(y)
 &\le C\ell^{-s},\\
\int_{(\Sigma\cap C_r)\setminus B_\ell(x)} K(x,y)\dd\sigma(y)
 &\le C\ell^{-1-s}.
\end{aligned}
\end{equation}

In the first line, for $|x-y|<\ell$ write
$|x-y|^{-m+1-s}=\ell^{-m+1-s}+(m-1+s)\int_{|x-y|}^{\ell}t^{-m-s}\dd t$.
Integrating and using \eqref{eq:all-growth} bounds the two terms by
$C\ell^{1-s}$ and $C\int_0^\ell t^{-s}\dd t$, respectively, before multiplication by $1-s$.
In the second line use
$|x-y|^{-m-s}=(m+s)\int_{|x-y|}^\infty t^{-m-s-1}\dd t$ on $|x-y|\ge\ell$. After interchanging the integrals, the area of $\Sigma\cap C_r$ with $\ell\le|x-y|\le t$ is at most $Ct^m$. This yields $C\int_\ell^\infty t^{-1-s}\dd t$.
For the last line integrate the kernel itself in the same way, with integrand $(m+1+s)t^{-m-2-s}$ in the representation, to get $C\int_\ell^\infty t^{-2-s}\dd t$. Thus the constants in all three estimates are uniform in $s$.

The next inequality is the truncated Jacobi inequality of Cabr\'e and Cozzi \cite[Theorem~1.3(ii)]{CC}. We use the normalization in \eqref{eq:operators} and track the constant in the cylinder truncation as $s\uparrow1$.

\begin{lemma}\label{lem:jacobi}
Let $E$ be a global subgraph, locally $s$-minimal in $C_{2r}$, with smooth boundary there. For $x\in\Sigma\cap C_{r/2}$,

\begin{equation}\label{eq:jacobi}
\begin{aligned}
L_{\Sigma\cap C_r}\omega(x)
&\ge (1-s)\omega(x)\int_{\Sigma\cap C_r}(1-\nu(x)\cdot\nu(y))K(x,y)\dd\sigma(y)\\
&\quad-C(m)(1-s)r^{-1-s}\omega(x).
\end{aligned}
\end{equation}
\end{lemma}
\begin{proof}

The curvature equation in the introduction gives $H_s[E](x)=0$ on the smooth minimizing boundary.

For a tangential vector $v$, Cabr\'e and Cozzi \cite[Proposition~4.1, (4.2)--(4.3)]{CC} differentiate along the boundary with $E$ fixed. The kernel derivative in the center is $-v\cdot D_yK(x,y)$. Their identity for a smooth boundary gives

\[
\partial_v H_s[E](x)
=-\PV\int_{\R^{m+1}}(\one_{E^c}-\one_E)(y)\,v\cdot D_yK(x,y)\dd y .
\]

For this principal value the cited proposition uses a radial cutoff and cancellation of the quadratic Taylor graph. In coordinates tangent at $x$, the quadratic graph is even in the base variables and the tangential kernel derivative is odd. The gap between the actual graph and its quadratic approximation is $O(|y-x|^{2+\beta})$ for any fixed $s<\beta<1$. It gives an error bounded by $C(E,s,\beta)\int_0^{2\delta}t^{\beta-s-1}\dd t$ in a cutoff at scale $\delta$. This tends to zero at fixed $s$.

Choose the tangential projection $v=e_{m+1}-\omega(x)\nu(x)$ of the vertical unit vector. Its components are $v'=-\omega(x)\nu'(x)$ and $v_{m+1}=1-\omega(x)^2$. Since $v\cdot\nu(x)=0$ and $H_s[E]=0$ along the graph, the preceding derivative is zero. Split the volume integral into the cylinder and its exterior. Inside the cylinder, $-\,(\one_{E^c}-\one_E)\operatorname{div}_y(Kv)$ is integrated over $E\cap C_r$ with a plus sign and $E^c\cap C_r$ with a minus sign. The surface fluxes there are

\[
\PV\int_{\Sigma\cap C_r} K v\cdot\nu(y)\dd\sigma(y)
-\PV\int_{\Sigma\cap C_r} K v\cdot(-\nu(y))\dd\sigma(y)
=2\PV\int_{\Sigma\cap C_r} K v\cdot\nu(y)\dd\sigma(y).
\]

The wall normal is $(y'/r,0)$. Splitting the exterior derivative into horizontal and vertical parts, as in \cite[(4.10)]{CC}, gives
\begin{equation}\label{eq:jacobi-identity}
\begin{aligned}
 2\PV\int_{\Sigma\cap C_r}(\omega(y)-\omega(x)\nu(x)\cdot\nu(y))K\dd\sigma(y)
 &=-\omega(x)\int_{\partial C_r}(\one_{E^c}-\one_E)\frac{\nu'(x)\cdot y'}r K\dd\HH^m\\
 &\quad-\omega(x)\int_{C_r^c}(\one_{E^c}-\one_E)\nu'(x)\cdot D_{y'}K\dd y\\
 &\quad+(1-\omega(x)^2)\int_{C_r^c}(\one_{E^c}-\one_E)\partial_{y_{m+1}}K\dd y .
\end{aligned}
\end{equation}
The trace on the wall is taken from inside. The cap fluxes vanish as the height tends to infinity. Since $|x'|<r/2$,
\[
\begin{aligned}
 \int_{\partial C_r}K\dd\HH^m
 &\le Cr^{m-1}\int_\R ((r/2)^2+t^2)^{-(m+1+s)/2}\dd t\le Cr^{-1-s},\\
 \int_{C_r^c}|D_{y'}K|\dd y
 &\le C\int_{r/2}^\infty t^{-2-s}\dd t\le Cr^{-1-s},\\
 \int_{C_r^c}(\one_{E^c}-\one_E)\partial_{y_{m+1}}K\dd y
 &=\int_{\R^m\setminus B'_r}\left(-\int_{-\infty}^{u(y')}+\int_{u(y')}^\infty\right)
       \partial_tK(x,(y',t))\dd t\dd y'\\
 &=-2\int_{\R^m\setminus B'_r}K(x,(y',u(y')))\dd y'\le0.
\end{aligned}
\]
Therefore \eqref{eq:jacobi-identity} yields
\[
 -\frac{L_{\Sigma\cap C_r}\omega(x)}{1-s}
 +\omega(x)\int_{\Sigma\cap C_r}(1-\nu(x)\cdot\nu(y))K\dd\sigma(y)
 \le C\omega(x)r^{-1-s}.
\]

\end{proof}

We derive a first-moment bound from the curvature equation $H_s[E]=0$ by comparing the flux on $\partial E$ with that on the boundary of its tangent half-space.

\begin{lemma}\label{lem:moment}
Fix a smooth nonincreasing radial cutoff $\eta$, equal to one on $B_{3/4}$ and zero outside $B_1$. Choose it as the square of a fixed smooth profile so that $|D\eta|\le C\sqrt\eta$ and $|D^2\eta|\le C$. We use the same profile on the base and in the ambient space. If $x\in\Sigma$ and $B_\rho(x)$ lies in the smooth minimizing region, then
\begin{equation}\label{eq:moment}
\begin{aligned}
&(1-s)\left|\PV\int_\Sigma \eta((y-x)/\rho)(y-x)K(x,y)\dd\sigma(y)\right|\\
&\quad\le C(1-s)\rho^{-s}
+C\rho(1-s)\int_\Sigma \eta((y-x)/\rho)(1-\nu(x)\cdot\nu(y))K(x,y)\dd\sigma(y).
\end{aligned}
\end{equation}
\end{lemma}

\begin{proof}

Translate and scale to $x=0$ and $\rho=1$. In this proof $r=|y|$, $K=K(0,y)$, and $H=\{y:y\cdot\nu(0)<0\}$ is the tangent half-space. Write $\dd\eta/\dd r$ for the derivative of the radial profile. The curvature equation and the symmetry of $H$ give $\int_{\R^{m+1}}(\one_E-\one_H)K\dd y=0$. The goal is to bound the vector moment $\PV\int_\Sigma\eta Ky\dd\sigma$. We use the projections of $y$ onto the tangent plane and its normal line.

Split the normal at the origin into its component along $y$ and its component orthogonal to $y$. In the normal difference this gives

\[
1-\nu(0)\cdot\nu(y)
=1-\left(\nu(0)-\frac{y\cdot\nu(0)}{r^2}y\right)\cdot\nu(y)
 -\frac{(y\cdot\nu(0))(y\cdot\nu(y))}{r^2}.
\]

After multiplication by $\eta K$, the last term is the negative boundary flux of $\eta K(y\cdot\nu(0))y/r^2$. Apply the divergence theorem to this field on $E$ and $H$ and subtract. The flux on $\partial H$ vanishes, and

\[
\Div\!\left(\eta K\frac{y\cdot\nu(0)}{r^2}y\right)
=-(1+s)\eta K\frac{y\cdot\nu(0)}{r^2}
 +\frac{\dd\eta}{\dd r}K\frac{y\cdot\nu(0)}{r}.
\]

Here the coefficient $-(1+s)$ results from $m$ in the divergence of $(y\cdot\nu(0))y/r^2$ and $-(m+1+s)$ from differentiating $K$. On $E\setminus H$ the indicator difference is one and $y\cdot\nu(0)\ge0$. On $H\setminus E$ both are negative. Thus $(\one_E-\one_H)(y\cdot\nu(0))=|y\cdot\nu(0)|\one_{E\triangle H}$. Substitution of the flux gives

\begin{equation}\label{eq:calibration}
\begin{aligned}
\int_\Sigma\eta K(1-\nu(0)\cdot\nu(y))\dd\sigma
={}&\int_\Sigma\eta K
\left[1-\left(\nu(0)-\frac{y\cdot\nu(0)}{r^2}y\right)\cdot\nu(y)\right]\dd\sigma\\
&+(1+s)\int_{E\triangle H}\eta K\frac{|y\cdot\nu(0)|}{r^2}\dd y\\
&-\int_{E\triangle H}\frac{\dd\eta}{\dd r}K\frac{|y\cdot\nu(0)|}{r}\dd y .
\end{aligned}
\end{equation}

The two volume terms are nonnegative since $\eta\ge0$ and $\dd\eta/\dd r\le0$. The orthogonal component of $\nu(0)$ has squared length $1-(y\cdot\nu(0))^2/r^2$, by the Pythagorean identity. Since $\nu(y)$ is a unit vector, its scalar product with that component is at most the length of the component. Subtracting this upper bound from one and rationalizing yields

\[
\begin{aligned}
1-\left(\nu(0)-\frac{y\cdot\nu(0)}{r^2}y\right)\cdot\nu(y)
&\ge 1-\sqrt{1-\frac{(y\cdot\nu(0))^2}{r^2}}\\
&=\frac{(y\cdot\nu(0))^2/r^2}{1+\sqrt{1-(y\cdot\nu(0))^2/r^2}}
\ge\frac{(y\cdot\nu(0))^2}{2r^2}.
\end{aligned}
\]

Thus the left side of \eqref{eq:calibration} controls half the integral of $\eta K(y\cdot\nu(0))^2/r^2$ over $\Sigma$, as well as the two nonnegative volume terms.

To estimate the moment in a unit tangent direction $e\perp\nu(0)$, consider the field $(e\cdot y)\nu(0)-(y\cdot\nu(0))e$. It has zero divergence and zero radial component. Its divergence remains zero when multiplied by the radial weight $\eta K$. Its flux on $\Sigma$ therefore gives

\[
\PV\int_\Sigma\eta K(e\cdot y)\dd\sigma
=\int_\Sigma\eta K\bigl[(e\cdot y)(1-\nu(0)\cdot\nu(y))
       +(y\cdot\nu(0))(e\cdot\nu(y))\bigr]\dd\sigma .
\]

The first term is bounded in absolute value by the left side of \eqref{eq:calibration} since $|e\cdot y|\le1$ on the support. By Cauchy--Schwarz the second is bounded by

\[
\left(\int_\Sigma\eta K\frac{(y\cdot\nu(0))^2}{r^2}\dd\sigma\right)^{1/2}
\left(\int_\Sigma\eta K r^2|e\cdot\nu(y)|^2\dd\sigma\right)^{1/2}.
\]

The first squared factor is controlled by the surface bound just obtained from \eqref{eq:calibration}. For the second use $r\le1$ and $|e\cdot\nu(y)|^2\le|\nu(y)-\nu(0)|^2=2(1-\nu(0)\cdot\nu(y))$. Thus both integrals under the square roots are bounded by twice the left side of \eqref{eq:calibration}.

To estimate the component along $\nu(0)$, apply the divergence theorem to $\eta K(y\cdot\nu(0))\nu(0)$. Differentiation in the fixed unit direction $\nu(0)$ gives $\partial_{\nu(0)}(y\cdot\nu(0))=1$ and $\partial_{\nu(0)}r=(y\cdot\nu(0))/r$. The latter follows by differentiating $|y+t\nu(0)|$ at $t=0$. The three product-rule terms are therefore

\[
\Div\bigl(\eta K(y\cdot\nu(0))\nu(0)\bigr)
=\eta K+K(y\cdot\nu(0))^2
  \left(\frac1r\frac{\dd\eta}{\dd r}-(m+1+s)\frac{\eta}{r^2}\right).
\]

The flux on $\Sigma$ contains the factor $\nu(0)\cdot\nu(y)$, and the flux on $\partial H$ vanishes. Apply the divergence theorem on $E$ minus that on $H$, and add the difference between the desired moment and the flux. This gives

\[
\begin{aligned}
\int_\Sigma\eta K(y\cdot\nu(0))\dd\sigma
={}&\int_\Sigma\eta K(y\cdot\nu(0))(1-\nu(0)\cdot\nu(y))\dd\sigma\\
&+\int_{\R^{m+1}}(\one_E-\one_H)\eta K\dd y\\
&+\int_{\R^{m+1}}(\one_E-\one_H)K(y\cdot\nu(0))^2
             \left(\frac1r\frac{\dd\eta}{\dd r}-(m+1+s)\frac{\eta}{r^2}\right)\dd y .
\end{aligned}
\]

The first term is bounded by the left side of \eqref{eq:calibration} since $|y\cdot\nu(0)|\le1$ on the support. By the curvature equation the second is $-\int(\one_E-\one_H)(1-\eta)K\dd y$, bounded by $|\mathbb S^m|\int_{3/4}^\infty t^{-1-s}\dd t$. In the last term, $|\one_E-\one_H|=\one_{E\triangle H}$ and $(y\cdot\nu(0))^2\le |y\cdot\nu(0)|$ on the support. Its absolute value is thus controlled by the two volume terms in \eqref{eq:calibration} with a factor at most $(m+1+s)/(1+s)$.

All divergence identities can first be computed with an inner radial cutoff at scale $\delta$. The tangential field has zero flux on the cutoff sphere. For the other two fields, the volume of $(E\triangle H)\cap B_{2\delta}$ is $O(\delta^{m+2})$ and $|y\cdot\nu(0)|=O(\delta^2)$ there. The cutoff derivatives cost $O_{E,s}(\delta^{1-s})$ and $O_{E,s}(\delta^{3-s})$, respectively. These tend to zero at fixed $s$. The tangential flux identity with a sphere deleted has an absolutely integrable right side, so its limit gives the principal value.

The tangential estimate holds for every unit $e$ in the tangent plane, so it bounds the norm of the tangential projection of the vector moment. The normal estimate bounds the other projection. At unit scale their sum is bounded by a constant plus a constant times the left side of \eqref{eq:calibration}. In the original coordinates the moment scales by $\rho^{-s}$ and the normal-difference integral by $\rho^{-1-s}$. This gives the factor $\rho$ in front of the latter integral. Multiplying by $1-s$ proves \eqref{eq:moment}.

\end{proof}

\section{Proof of the main theorem}\label{sec:proof}
Both estimates in the proof concern the same graph patch $D=\Sigma\cap B_{4r}(P)$ and the same function $g=\omega^{-q}-1$. Lemma~\ref{lem:radial-mean} bounds $g(P)$ by $r^{-m}\int_D g\,\dd\sigma$ and a truncation error. Lemma~\ref{lem:height} bounds that integral average by $C(m)q(1+\osc_{B'_{20r}(P')}u/r)$. We prove these two local estimates first.

\subsection{The radial potential and the pointwise mean-value inequality}

Translate $P\in\Sigma$ to the origin and take the outer cutoff radius to be one. Suppose $E$ is $s$-minimal in $B_{40}$ with smooth boundary there. Set $D=\Sigma\cap B_4$. Lemma~\ref{lem:geometry} gives its area growth bound. Indeed, for a small ball meeting $D$, recenter at a point of the intersection and apply \eqref{eq:area} to a ball of twice the radius. For larger radii use the total area of $D$. The radial-layer argument for \eqref{eq:moments-upper} thus applies on $D$. We first compute the planar potential, then compare the surface integral with the tangent-plane integral, and finally extract the value of the test function at the pole.

Let $T$ be an $m$-plane through the origin. Use the fundamental solution of $(-\Delta)^{(1+s)/2}$ on $T$ from \cite[Section~2]{Bucur}. Multiplying its coefficient by the usual fractional Laplacian kernel constant and dividing by $1-s$ gives the normalization for \eqref{eq:operators}. Thus $G_s(x)=c_{m,s}|x|^{1+s-m}$ satisfies
\[
\int_T G_s\,L_T\varphi\,\dd\HH^m=\varphi(0)
\qquad\text{for every }\varphi\in C_c^\infty(T).
\]
With the cutoff of Lemma~\ref{lem:moment}, define
\begin{equation}\label{eq:radial-phi}
\Phi(x)=c_{m,s}\eta(x)^2\bigl(|x|^{1+s-m}-2^{1+s-m}\bigr).
\end{equation}
Near the pole $\Phi-G_s$ is constant.

Using $\Gamma(z+1)=z\Gamma(z)$ in the coefficient just described gives
\[
(m-1-s)c_{m,s}
=\frac{(1+s)\Gamma((m+1+s)/2)\Gamma((m+1-s)/2)}
 {2\pi^m\Gamma((3-s)/2)\Gamma((1+s)/2)} .
\]
The product $(m-1-s)c_{m,s}$ has a positive lower bound and a finite upper bound depending only on $m$. Indeed, all the gamma arguments stay in compact subsets of $(0,\infty)$ for fixed $m\ge2$ and $3/4\le s\le1$. For $m=2$ specifically,
$(1-s)c_{2,s}=(1+s)^2/(4\pi^2)$.
Consequently, in dimension two, $\Phi(x)$ tends to $\pi^{-2}\eta(x)^2\log(2/|x|)$ for $x\ne0$ as $s\uparrow1$. For every $m\ge2$ and $0<|x|\le1/2$,
\begin{equation}\label{eq:radial-phi-basic}
\Phi(x)=(m-1-s)c_{m,s}\int_{|x|}^2 t^{s-m}\dd t,\qquad
\Phi(x)\ge c |x|^{1+s-m},\qquad \int_D\Phi\dd\sigma\le C .
\end{equation}
For the lower bound integrate over $|x|<t<2|x|$. For the mass, Tonelli gives
\[
\int_D\Phi\dd\sigma
 \le C\int_0^2 t^{s-m}\HH^m(D\cap B_t)\dd t
 \le C\int_0^2 t^s\dd t\le C .
\]

The fundamental-solution identity gives $L_TG_s=0$ off the origin on a plane through the origin. For an affine plane at positive distance from the origin, the calculation uses stereographic projection.

\begin{lemma}\label{lem:offset}
For every affine $m$-plane $T\subset\R^{m+1}$ and $x\in T\setminus\{0\}$,
\[
 L_TG_s(x)=b_{m,s}\dist(0,T)^{1+s}|x|^{-m-1-s},\qquad b_{m,s}>0 .
\]
\end{lemma}

\begin{proof}
If $0\in T$, scale to $|x|=1$. Split the integral at $|y|=1$ and invert the inner region. For $0<\delta<1/2$,
\[
\begin{aligned}
&\int_{T\setminus B_\delta(x)}(1-|y|^{1+s-m})K(x,y)\dd y\\
&=\int_{T\cap\{|y|>1\}} (1-|y|^{1+s-m})K(x,y)
 \bigl(\one_{\{|x-y|>\delta\}}-\one_{\{|x-y|>\delta |y|\}}\bigr)\dd y .
\end{aligned}
\]

The indicators differ only for $\delta<|x-y|\le\delta|y|\le\delta(1+|x-y|)$, hence $|x-y|\le\delta/(1-\delta)$. Its thickness is $\delta/(1-\delta)-\delta=\delta^2/(1-\delta)\le2\delta^2$. The shell has $m$-volume $O(\delta^{m+1})$, and $|1-|y|^{1+s-m}|\le C|x-y|$ on it. Its integral is $O(\delta^{1-s})$. Letting $\delta\downarrow0$ gives $L_TG_s(x)=0$. If $0\notin T$, rotate and scale to $T=\{(z,1):z\in\R^m\}$. In its base coordinates the inverse stereographic projection is $\zeta(z)=(2z,|z|^2-1)/(1+|z|^2)$. In particular,
\[
|\zeta(x)-\zeta(y)|=\frac{2|x-y|}{(1+|x|^2)^{1/2}(1+|y|^2)^{1/2}},\qquad
\dd\HH^m(\zeta(y))=\frac{2^m}{(1+|y|^2)^m}\dd y .
\]
It follows that
\[
\int_{\R^m}\frac{|x-y|^{1+s-m}}{(1+|y|^2)^{(m+1+s)/2}}\dd y
=b_{m,s}^{-1}(1+|x|^2)^{-(m-1-s)/2}.
\]
The angular integral $\int_{\mathbb S^m}|\zeta(x)-\theta|^{1+s-m}\dd\HH^m(\theta)$ is finite and independent of $x$, hence $b_{m,s}>0$. Applying the planar operator to the convolution in the sense of distributions and using $L_{\R^m}G_s=\delta_0$ gives
\[
 L_{\R^m}\!\left[c_{m,s}(1+|\cdot|^2)^{-(m-1-s)/2}\right](x)
 = b_{m,s}(1+|x|^2)^{-(m+1+s)/2}.
\]
Both sides are smooth. Dilation to height $h=\dist(0,T)>0$ multiplies the operator value by $h^{-m}$ and replaces the base coordinate $x$ by $x/h$. This proves the formula.
\end{proof}

We use Lemma~\ref{lem:offset} to compare the operator applied to $\Phi$ on the surface with the planar operator and estimate the error of this comparison.

\begin{lemma}\label{lem:defect}
Let $0<q\le1/2$ and $x\in D\setminus\{0\}$. All constants below are dimensional. If $0<|x|<1/4$,
\begin{equation}\label{eq:bad-defect}
\begin{aligned}
&(1-s)\PV\int_D \frac{1+\nu(x)\cdot\nu(y)}{2}\bigl(\Phi(x)-\Phi(y)\bigr)K(x,y)\dd\sigma(y)\\
&\quad\ge-\frac{q(1-s)}4\,\Phi(x)\int_D (1-\nu(x)\cdot\nu(y))K(x,y)\dd\sigma(y)-C(1-s)q^{-s}|x|^{-m}.
\end{aligned}
\end{equation}
If $E$ satisfies \eqref{eq:flatness} in $B_\ell$ with $0<\ell\le1/4$, the last term can be replaced for $0<|x|<c_1\ell$ by
\begin{equation}\label{eq:good-defect}
-C(1-s)q^{-s}|x|^{-m}(|x|/\ell)^{1/[2(m+1)]}.
\end{equation}
Here $c_1$ is dimensional. For $|x|\ge1/4$ and $(1-s)q^{-s}\le1$ one has
\begin{equation}\label{eq:outer-defect}
\begin{aligned}
&(1-s)\PV\int_D \frac{1+\nu(x)\cdot\nu(y)}{2}\bigl(\Phi(x)-\Phi(y)\bigr)K(x,y)\dd\sigma(y)\\
&\quad\ge-\frac{q(1-s)}4\,\Phi(x)\int_D (1-\nu(x)\cdot\nu(y))K(x,y)\dd\sigma(y)-C.
\end{aligned}
\end{equation}
\end{lemma}

\begin{proof}
Let $T_x$ be the tangent plane to $\Sigma$ at $x$ and set $H_x=\{y:(y-x)\cdot\nu(x)<0\}$. For $\varphi\in C_c^1(B_4)$, applying the divergence theorem to $E$ and $H_x$ gives
\[
 \int_E\partial_{\nu(x)}\varphi\dd y
 =\int_\Sigma\varphi\,\nu(x)\cdot\nu(y)\dd\sigma(y),\qquad
 \int_{H_x}\partial_{\nu(x)}\varphi\dd y=\int_{T_x}\varphi\dd\HH^m .
\]
Subtract and add $\frac12\int_\Sigma(1-\nu(x)\cdot\nu(y))\varphi\dd\sigma(y)$ to obtain
\begin{equation}\label{eq:tangent-diff}
\begin{aligned}
&\int_\Sigma \frac{1+\nu(x)\cdot\nu(y)}2 \varphi\dd\sigma(y)-\int_{T_x}\varphi\dd\HH^m\\
&\quad=\frac12\int_\Sigma (1-\nu(x)\cdot\nu(y))\varphi\dd\sigma(y)
+\int_{\R^{m+1}}(\one_E-\one_{H_x})\,\partial_{\nu(x)}\varphi\dd y .
\end{aligned}
\end{equation}
The curvature equation gives $\int_{\R^{m+1}}(\one_E-\one_{H_x})K\dd y=0$. Inserting the singular tests below is justified by radial cutoffs. In $B_{2\delta}(x)\cap(E\triangle H_x)$ smooth tangency gives $|(y-x)\cdot\nu(x)|\le C_x\delta^2$ and volume $\le C_x\delta^{m+2}$. For a cutoff at $x$ its normal derivative is $O(|(y-x)\cdot\nu(x)|/\delta^2)$ on the annulus, so the error for a potential difference times $K$ is $O_x(\delta^{2-s})$. At the pole $0\ne x$ the cutoff error is $O_{x,s}(\delta^{1+s})$.

In \eqref{eq:calibration} use the cutoff $\eta((y-x)/(2\rho))$, with $B_{2\rho}(x)\subset B_4$. It equals one on $B_\rho(x)$ and all three terms on the right are nonnegative. Hence
\begin{equation}\label{eq:body-moment}
\begin{aligned}
 &(1+s)(1-s)\int_{B_\rho(x)\cap(E\triangle H_x)}
        \frac{|(y-x)\cdot\nu(x)|}{|x-y|^2}\,K(x,y)\dd y\\
 &\quad\le (1-s)\int_D\eta((y-x)/(2\rho))(1-\nu(x)\cdot\nu(y))K(x,y)\dd\sigma(y)\\
 &\quad\le (1-s)\int_D(1-\nu(x)\cdot\nu(y))K(x,y)\dd\sigma(y).
\end{aligned}
\end{equation}
The corresponding first absolute moment satisfies
\begin{equation}\label{eq:body-moment-linear}
(1-s)\int_{B_\rho(x)\cap(E\triangle H_x)}|x-y|\,K\dd y
\le C(1-s)\rho^{(3-s)/2}\left[\int_D (1-\nu(x)\cdot\nu(y))K(x,y)\dd\sigma(y)\right]^{1/2}.
\end{equation}
To prove this, set $U=(E\triangle H_x)\cap\{t/2<|x-y|<t\}$. Slicing parallel to $T_x$ gives, for $h>0$,
\[
 |U\cap\{|(y-x)\cdot\nu(x)|\le h\}|\le\int_{-h}^h |B'_t|\dd v=2|B'_1|t^m h .
\]

With $h=|U|/(4|B'_1|t^m)$ one has $\int_U|(y-x)\cdot\nu(x)|\dd y\ge h|U|/2$. Using $t/2<|x-y|<t$ gives
\[
(1-s)\int_U |x-y|K\dd y\le C(1-s)t^{-m-s}|U|\le C(1-s)t^{(3-s)/2}
 \left(\int_U\frac{|(y-x)\cdot\nu(x)|}{|x-y|^2}\,K\dd y\right)^{1/2}.
\]
Sum over $t=2^{-j}\rho$, $j\ge0$. Cauchy--Schwarz, $\sum_{j\ge0}(2^{-j}\rho)^{3-s}\le C\rho^{3-s}$ and \eqref{eq:body-moment} give \eqref{eq:body-moment-linear}.

Write $r=|x|$. For $0<r\le1/4$, set $\rho=cqr$, with $c>0$ sufficiently small depending only on $m$. We first estimate the integral with the cutoff $\eta((y-x)/\rho)$. On the support of $\eta((y-x)/\rho)$, the difference $\Phi(x)-\Phi(y)$ equals $G_s(x)-G_s(y)$.
Here $|z|\asymp r$ on the segment from $x$ to $y$, so $|G_s(x)-G_s(y)|\le C r^{s-m}|x-y|$. In \eqref{eq:tangent-diff} insert $\varphi(y)=\eta((y-x)/\rho)(G_s(x)-G_s(y))K(x,y)$. With $\partial_{\nu(x)}K=-(m+1+s)((y-x)\cdot\nu(x))K/|x-y|^2$, \eqref{eq:body-moment} gives
\[
\begin{aligned}
 &\tfrac12\int_\Sigma (1-\nu(x)\cdot\nu(y))|\varphi(y)|\dd\sigma(y)\\
 &\quad+\int_{E\triangle H_x}\eta((y-x)/\rho)|G_s(x)-G_s(y)|\,|\partial_{\nu(x)}K|\dd y\\
 &\le Cr^{s-m}\rho\int_D(1-\nu(x)\cdot\nu(y))K\dd\sigma(y),\\
 &\int_{E\triangle H_x}|\partial_{\nu(x)}[\eta((y-x)/\rho)]|\,|G_s(x)-G_s(y)|K\dd y\\
 &\le Cr^{s-m}\int_{\rho/2<|x-y|<\rho}K\dd y\le Cr^{s-m}\rho^{-s}.
\end{aligned}
\]

Apart from the factor $(m-1-s)c_{m,s}$, differentiation of the potential gives
\[
(1-s)\int(\one_E-\one_{H_x})\eta((y-x)/\rho) |y|^{s-m-1}\bigl(\nu(x)\cdot x+(y-x)\cdot\nu(x)\bigr)K\dd y.
\]
The contribution of $(y-x)\cdot\nu(x)$ is nonnegative since $(\one_E-\one_{H_x})((y-x)\cdot\nu(x))\ge0$. In the part with $\nu(x)\cdot x$, subtract $r^{s-m-1}$ from $|y|^{s-m-1}$. For the constant part we use
\[
\begin{aligned}
&\int(\one_E-\one_{H_x})\eta((y-x)/\rho) K\dd y\\
&\quad=-\int(\one_E-\one_{H_x})(1-\eta((y-x)/\rho))K\dd y,
\end{aligned}
\]
whose right side is bounded in absolute value by $C\rho^{-s}$. With $|\nu(x)\cdot x|\le r$ this gives $C(1-s)r^{s-m}\rho^{-s}$. The difference is bounded by $Cr^{s-m-2}|x-y|$. By \eqref{eq:body-moment-linear} and Young's inequality it gives
\[
\begin{aligned}
&C|\nu(x)\cdot x|r^{s-m-2}(1-s)\rho^{(3-s)/2}\left[\int_D (1-\nu(x)\cdot\nu(y))K(x,y)\dd\sigma(y)\right]^{1/2}\\
&\quad\le \frac{q(1-s)}8\,\Phi(x)\int_D (1-\nu(x)\cdot\nu(y))K(x,y)\dd\sigma(y)\\
&\qquad+C(m)(1-s)q^{-1}r^{-m}\left(\frac{\nu(x)\cdot x}{r}\right)^2(\rho/r)^{3-s}.
\end{aligned}
\]
The lower bound in \eqref{eq:radial-phi-basic} is used in the Young estimate. Take $\rho/(qr)$ small enough that the surface term and the kernel-derivative term together are bounded by the same $q/8$ multiple of the normal-difference term. This gives the first term on the right of \eqref{eq:bad-defect}. Substituting $\rho=cqr$ bounds the remaining terms by $C(1-s)q^{-s}r^{-m}$.

When \eqref{eq:flatness} holds in $B_\ell$ and $r<c_1\ell$, the common graph in $B_{\vartheta\ell/4}(x)$ satisfies
\begin{equation}\label{eq:local-chart}
\begin{gathered}
|\nu(y)-\nu(x)|\le C(|x-y|/\ell)^\gamma,\qquad
|\nu(x)\cdot x|\le C r(r/\ell)^\gamma,\\
|(y-x)\cdot\nu(x)|\le C |x-y|(|x-y|/\ell)^\gamma\qquad\bigl(y\in E\triangle H_x\bigr).
\end{gathered}
\end{equation}
These bounds follow from Lemma~\ref{lem:flat}. The normal difference is evaluated at $y\in\Sigma$. Here $\vartheta$ is the radius factor in that lemma, and we take $c_1<\vartheta/400$. From the volume on shells between $\rho/2$ and $\vartheta\ell/4$,
\[
\int_{E\triangle H_x}(1-\eta((y-x)/\rho))K\dd y\le C\rho^{-s}(\rho/\ell)^\gamma .
\]
The part outside $B_{\vartheta\ell/4}(x)$ is bounded by $C\ell^{-s}$. The cutoff-derivative term gains the same factor as the displayed tail. The nonpotential terms therefore gain $(r/\ell)^\gamma$, with the Young remainder gaining $(r/\ell)^{2\gamma}$ from $\nu(x)\cdot x$.

We next estimate the integral with the complementary cutoff $1-\eta((y-x)/\rho)$.
In the general case, bound the negative part of the surface contribution for $|x-y|>\rho/2$ and the absolute value of the corresponding full planar tail. For $|x-y|<r/2$ the difference is at most $Cr^{s-m}|x-y|$, giving a bound $Cr^{s-m}\rho^{-s}$ by area growth. For $|x-y|\ge r/2$ and $|y|\le4r$ use
\begin{equation}\label{eq:integrable-riesz}
\int_{\Sigma\cap B_{5r}}\big|G_s(x)-G_s(y)\big|\dd\sigma(y)\le Cr^{1+s}.
\end{equation}
The same bound holds for planar area. Write the absolute difference as $(m-1-s)c_{m,s}$ times the integral of $t^{s-m}$ between the two radii. For $|y|\le r$ use Tonelli and growth, and for $r\le|y|\le5r$ use the integral over these comparable radii. For more distant points the surface contribution is nonnegative by radial monotonicity of $\Phi$. Use $|G_s(x)-G_s(y)|\le Cr^{1+s-m}(1+\log(|x-y|/r))$ on that planar tail, of size $Cr^{-m}$. Adding the integrals with cutoffs $\eta((y-x)/\rho)$ and $1-\eta((y-x)/\rho)$ gives
\[
\begin{aligned}
 &(1-s)\PV\int_D\tfrac{1+\nu(x)\cdot\nu(y)}2(\Phi(x)-\Phi(y))K\dd\sigma(y)\\
 &\quad+\tfrac q4(1-s)\Phi(x)\int_D(1-\nu(x)\cdot\nu(y))K\dd\sigma(y)\\
 &\ge L_{T_x}G_s(x)-C(1-s)(r^{s-m}\rho^{-s}+r^{-m}+q^{-1}r^{-m}(\rho/r)^{3-s})\\
 &\ge -C(1-s)q^{-s}r^{-m},
\end{aligned}
\]
since $L_{T_x}G_s(x)\ge0$ by Lemma~\ref{lem:offset} and $\rho=cqr$. This proves \eqref{eq:bad-defect}.

When \eqref{eq:flatness} holds in $B_\ell$, use the chart cutoff $\eta(4(y-x)/(\vartheta\ell))$ in the transfer formula \eqref{eq:tangent-diff}, with
$(\eta(4(y-x)/(\vartheta\ell))-\eta((y-x)/\rho))(G_s(x)-G_s(y))K$.

For $\rho/2<|x-y|<r/2$ the cut-off integrand is bounded by $Cr^{s-m}|x-y|K$ and its derivative by $Cr^{s-m}K$. In the surface term use the squared normal difference. In the volume term a shell of radius $t$ has relevant volume at most $C(t/\ell)^\gamma t^{m+1}$. Their sum is bounded by $Cr^{s-m}\rho^{-s}(r/\ell)^\gamma$.

For $r/2\le |x-y|\le4r$ the two cutoffs are constant, with difference one. The surface term is bounded by $Cr^{-m}(r/\ell)^{2\gamma}$ using \eqref{eq:integrable-riesz}. In the volume term $|D_y[(G_s(x)-G_s(y))K]|\le Cr^{-m-1-s}|y|^{s-m}$ and its support has volume at most $Cr^{m+1}(r/\ell)^\gamma$. The bound $\int_U|y|^{s-m}\dd y\le C|U|^{(1+s)/(m+1)}$ follows by splitting at radius $|U|^{1/(m+1)}$. This gives $Cr^{-m}(r/\ell)^{\gamma(1+s)/(m+1)}$.

For $4r<|x-y|<\vartheta\ell/4$, $|y|\asymp |x-y|$ and $|G_s(x)-G_s(y)|\le Cr^{1+s-m}(1+\log(|x-y|/r))$, $|DG_s(y)|\le C|x-y|^{s-m}$. On a shell of radius $t$, the surface estimate is $C(t/\ell)^{2\gamma}r^{1+s-m}t^{-1-s}(1+\log(t/r))$. The volume estimate has $\gamma$ in place of $2\gamma$. Since $2\gamma<1+s$, their sum is at most $Cr^{-m}(r/\ell)^\gamma$. Outside the chart cutoff the surface term is nonnegative and the planar tail is bounded by $Cr^{1+s-m}\ell^{-1-s}(1+\log(\ell/r))$. Including $1-s$ and the near-diagonal bound proves \eqref{eq:good-defect}.

For $r=|x|\ge1/4$, the square cutoff in \eqref{eq:radial-phi} gives $|D\Phi(x)|\le C\sqrt{\Phi(x)}$, and $|D^2\Phi|\le C$ on $\{|y|\ge1/8\}$. If $\Phi(x)>0$, take $\rho=cq\sqrt{\Phi(x)}<1/32$. In the linear Taylor term write $(1+\nu(x)\cdot\nu(y))/2=1-(1-\nu(x)\cdot\nu(y))/2$. By \eqref{eq:moment},
\[
\begin{aligned}
 &(1-s)\left|\PV\int_D\eta((y-x)/\rho)\tfrac{1+\nu(x)\cdot\nu(y)}2
        D\Phi(x)\cdot(y-x)K(x,y)\dd\sigma(y)\right|\\
 &\quad\le C(1-s)|D\Phi(x)|\rho^{-s}
       +C(1-s)\rho|D\Phi(x)|\int_D(1-\nu(x)\cdot\nu(y))K(x,y)\dd\sigma(y).
\end{aligned}
\]
Choose $c=c(m)$ so that $C\rho |D\Phi(x)|\le q\Phi(x)/4$. The quadratic remainder on $B_{1/16}(x)$, the negative contribution outside it, and the linear term where $1-\eta((y-x)/\rho)$ is nonzero are controlled by \eqref{eq:moments-upper} and \eqref{eq:radial-phi-basic} as follows.
\[
\begin{aligned}
 &(1-s)\int_{D\cap B_{1/16}(x)}|x-y|^2K\dd\sigma(y)
 +(1-s)\int_{D\setminus B_{1/16}(x)}\Phi(y)K\dd\sigma(y)\le C,\\
 &(1-s)|D\Phi(x)|\int_{D\setminus B_{\rho/2}(x)}|x-y|K\dd\sigma(y)
 \le C(1-s)\sqrt{\Phi(x)}\,\rho^{-s}\le C(1-s)q^{-s}.
\end{aligned}
\]
If $\Phi(x)=0$ then $D\Phi(x)=0$, so the first line alone suffices. This proves \eqref{eq:outer-defect}.
\end{proof}

\noindent\begin{minipage}{\linewidth}
Lemma~\ref{lem:defect} estimates the operator applied to $\Phi$ away from the pole. The following identity also accounts for the contribution at the pole.

\begin{lemma}\label{lem:dirac}
Fix the smooth surface and $s$. The principal value on the left of \eqref{eq:bad-defect}, considered as a function on $D\setminus\{0\}$, is in $L^1(D,\dd\sigma)$. For $f$ smooth on a surface neighborhood of $\overline D$,
\begin{equation}\label{eq:radial-dirac}
\begin{aligned}
&(1-s)\int_D \Phi(x)\PV\int_D \frac{1+\nu(x)\cdot\nu(y)}{2}\bigl(f(x)-f(y)\bigr)K(x,y)\dd\sigma(y)\dd\sigma(x)\\
&\quad= f(0)+(1-s)\int_D f(x)\PV\int_D \frac{1+\nu(x)\cdot\nu(y)}{2}\bigl(\Phi(x)-\Phi(y)\bigr)K(x,y)\dd\sigma(y)\dd\sigma(x).
\end{aligned}
\end{equation}
\end{lemma}
\end{minipage}

\begin{proof}
Fix $s$ and the surface. Write $I(f)$ for the integral on the left of \eqref{eq:radial-dirac}, and set, for $x\ne0$,
\[
p_\Phi(x)=(1-s)\PV\int_D\tfrac{1+\nu(x)\cdot\nu(y)}2
              (\Phi(x)-\Phi(y))K(x,y)\dd\sigma(y).
\]
Smoothness at $0$ gives a fixed $0<\ell\le1/4$ for which \eqref{eq:flatness} holds in $B_\ell$. Near $0$ the integral of $(1-\nu(x)\cdot\nu(y))K(x,y)$ over $D$ is bounded at fixed $s$, since $1-\nu(x)\cdot\nu(y)=\tfrac12|\nu(x)-\nu(y)|^2\le C_\Sigma |x-y|^2$ locally. Taking $q=1/2$ in \eqref{eq:bad-defect}--\eqref{eq:good-defect} and using $\Phi(x)\le c_{m,s}|x|^{1+s-m}$ gives, for $0<b<c_1\ell$,
\[
\begin{aligned}
 p_\Phi^-(x)&\le C_{\Sigma,s}\bigl(|x|^{1+s-m}+|x|^{-m+1/[2(m+1)]}\bigr)
      &&(0<|x|<c_1\ell),\\
 \int_{D\cap B_b}p_\Phi^-\,\dd\sigma
 &\le C_{\Sigma,s}\bigl(b^{1+s}+b^{1/[2(m+1)]}\bigr),
\end{aligned}
\]
where $p_\Phi^\pm=\max\{\pm p_\Phi,0\}$ and the second estimate follows by summing over dyadic annuli. Away from $0$ the principal value is locally bounded by Taylor expansion. Near $\partial D$, the support of $\Phi$ is at positive distance. Thus $p_\Phi^-\in L^1(D)$.
If a smooth test $\psi$ vanishes near $0$, exchange $x,y$ on $|x-y|>\epsilon$ and let $\epsilon\downarrow0$. Taylor expansion on $\operatorname{supp}\Phi$ bounds the truncated operator applied to $\psi$, and $\int_D\Phi\,\dd\sigma<\infty$ controls the pole in the exchanged integral. Thus
\[
I(\psi)=\int_D\psi p_\Phi\,\dd\sigma .
\]
To compute the contribution at $0$, put $D_t=t^{-1}D$ and let $\dd\sigma_t$ be its area measure. In $I(f\eta(\cdot/t))$ the changes $x=tX$, $y=tY$ give
\[
\begin{aligned}
I(f\eta(\cdot/t))={}&(1-s)\int_{D_t}t^{m-1-s}\Phi(tX)\\[-2pt]
&{}\cdot\PV\int_{D_t}\tfrac{1+\nu(tX)\cdot\nu(tY)}2
 \frac{f(tX)\eta(X)-f(tY)\eta(Y)}{|X-Y|^{m+1+s}}
 \dd\sigma_t(Y)\dd\sigma_t(X).
\end{aligned}
\]
Over the tangent plane $T_0$ write $\Sigma$ as the graph of $v$, with $v(0)=Dv(0)=0$. On each fixed disk $\|v(t\,\cdot)/t\|_{C^2(B'_b)}\le C_b t$. Hence the area densities and normal weights tend to one, $f(t\,\cdot)\eta\to f(0)\eta$ in $C^2$ in these charts, and $t^{m-1-s}\Phi(tX)\to G_s(X)$ off $0$. On $B_b$ the error from deleting $|X-Y|<\epsilon$ in the inner principal value is at most $C_{b,\Sigma,s}\int_0^\epsilon \rho^{-s}\dd\rho$. Indeed, in coordinates tangent at $X$ the linear Taylor term is odd and the remainder is quadratic. The area bound $\sigma_t(B_a(Z))\le Ca^m$ gives, uniformly in $t$, for $0<\delta<1$ and $M>2b>4$,
\[
\begin{aligned}
 \int_{D_t\cap B_\delta}|X|^{1+s-m}\dd\sigma_t(X)&\le C_s\delta^{1+s},\\
 \sup_{|X|\le b}\int_{D_t\setminus B_M}|X-Y|^{-m-1-s}\dd\sigma_t(Y)&\le C_sM^{-1-s},\\
 \int_{D_t\setminus B_b}|X|^{1+s-m}
       \int_{D_t\cap B_1}K(X,Y)\dd\sigma_t(Y)\dd\sigma_t(X)
 &\le C\int_{D_t\setminus B_b}|X|^{-2m}\dd\sigma_t(X)\le Cb^{-m}.
\end{aligned}
\]
Since $0\le t^{m-1-s}\Phi(tX)\le c_{m,s}|X|^{1+s-m}$, these bounds justify passage to the planar integral. The normalization of $G_s$ gives
\[
\lim_{t\downarrow0} I(f\eta(\cdot/t))=f(0)\int_{T_0}G_s L_{T_0}\eta\,\dd\HH^m=f(0).
\]
Take $f=1$. Since $I(1)=0$, the identity for tests vanishing near $0$ gives $\int_D(1-\eta(\cdot/t))p_\Phi\,\dd\sigma=-I(\eta(\cdot/t))\to-1$. Fatou's lemma and the integrability of $p_\Phi^-$ give
\[
 \int_D p_\Phi^+\dd\sigma\le\liminf_{t\downarrow0}\int_D(1-\eta(\cdot/t))p_\Phi^+\dd\sigma
 =-1+\int_D p_\Phi^-\dd\sigma<\infty .
\]
Thus $p_\Phi\in L^1(D)$. For general $f$ we may now take the limit in
\[
I(f)=I(f\eta(\cdot/t))+\int_D f(1-\eta(\cdot/t))p_\Phi\,\dd\sigma
 \ \longrightarrow\ f(0)+\int_D f p_\Phi\,\dd\sigma .
\]
\end{proof}

\noindent\begin{minipage}{\linewidth}
The sign of the planar term from Lemma~\ref{lem:offset} and the comparison in Lemma~\ref{lem:defect} control the surface error. The identity in Lemma~\ref{lem:dirac} gives the point value when we transfer the operator to $\Phi$. We combine these facts with the truncated Jacobi inequality to obtain the pointwise estimate for $g$.

\begin{lemma}\label{lem:radial-mean}
Let $s$ be above the threshold of Lemma~\ref{lem:nonflat}. Suppose the global subgraph $E$ is locally $s$-minimal in $C_{40r}(P')$ and its graph is smooth there. Let $P=(P',u(P'))$, $D=\Sigma\cap B_{4r}(P)$ and $\omega=\nu_{m+1}$. If $g=\omega^{-q}-1$, $0<q\le1/2$, $0\le g\le e^2-1$ on $D$, and $(1-s)q^{-1-s}\le c(m)$, then
\begin{equation}\label{eq:radial-mean}
g(P)\le C(m)\left[r^{-m}\int_D g\dd\sigma+(1-s)q^{-s}\right].
\end{equation}

\end{lemma}
\end{minipage}

\begin{proof}
Translate and rescale to $P=0$, $r=1$. For $x\in\Sigma\cap B_1$, restrict \eqref{eq:jacobi} from $\Sigma\cap C_4$ to $D=\Sigma\cap B_4$. Because $|x-y|\ge3$ on the discarded set, \eqref{eq:moments-upper} gives
\[
\begin{aligned}
 L_D\omega(x)&-(1-s)\omega(x)\int_D(1-\nu(x)\cdot\nu(y))K\dd\sigma(y)\\
 &\ge -C(1-s)\omega(x)\\
 &\quad-(1-s)\int_{(\Sigma\cap C_4)\setminus D}
     [\omega(x)\nu(x)\cdot\nu(y)-\omega(y)]K\dd\sigma(y)
 \ge -C(1-s)\omega(x).
\end{aligned}
\]
We next use this bound for $\omega$ to derive \eqref{eq:sym-jacobi} for $g=\omega^{-q}-1$.
Fix $x$ in $\Sigma\cap B_1$ and $y\in D$. Write $t=\omega(y)/\omega(x)>0$ for this calculation. Since $1+g(x)=\omega(x)^{-q}$,
\[
 g(x)-g(y)=(1+g(x))(1-t^{-q}).
\]
The second derivative of $t^{-q}$ is $q(q+1)t^{-q-2}>0$. Its tangent at $1$ therefore gives $t^{-q}\ge1-q(t-1)$. We multiply the resulting upper bound on $g(x)-g(y)$ by the nonnegative symmetric weight $(1+\nu(x)\cdot\nu(y))/2$. Since $t>0$ and $\nu(x)\cdot\nu(y)\le1$,
\[
 \frac{1+\nu(x)\cdot\nu(y)}2(t-1)
 = (t-1)+\frac{1-\nu(x)\cdot\nu(y)}2(1-t)
 \le (t-1)+\frac{1-\nu(x)\cdot\nu(y)}2 .
\]
By the identity $1-\nu(x)\cdot\nu(y)=|\nu(x)-\nu(y)|^2/2$ the normal difference is integrable against $K$ near $x$. The bound just proved for $L_D\omega$ gives
\[
\begin{aligned}
 &(1-s)\PV\int_D\left(t-1+\tfrac{1-\nu(x)\cdot\nu(y)}2\right)K(x,y)\dd\sigma(y)\\
 &= -\frac{L_D\omega(x)}{\omega(x)}
    +\frac{1-s}{2}\int_D(1-\nu(x)\cdot\nu(y))K(x,y)\dd\sigma(y)\\
 &\le C(m)(1-s)-\frac{1-s}{2}\int_D(1-\nu(x)\cdot\nu(y))K(x,y)\dd\sigma(y).
\end{aligned}
\]
Multiplying by $q(1+g(x))$ and using $g(x)-g(y)\le q(1+g(x))(t-1)$ yields
\begin{equation}\label{eq:sym-jacobi}
\begin{aligned}
&(1-s)\PV\int_D \frac{1+\nu(x)\cdot\nu(y)}{2}\bigl(g(x)-g(y)\bigr)K(x,y)\dd\sigma(y)\\
&\quad+\frac{q(1-s)}2(1+g(x))\int_D (1-\nu(x)\cdot\nu(y))K(x,y)\dd\sigma(y)\\
&\quad\le C(m)q(1-s)(1+g(x)).
\end{aligned}
\end{equation}

We use the radial potential $\Phi$ from \eqref{eq:radial-phi} to control the logarithmic term that arises when we integrate \eqref{eq:bad-defect}--\eqref{eq:outer-defect} against $g$. Choose $0<\ell\le1/(16A)$ such that $E$ satisfies \eqref{eq:flatness} in $B_\ell$ and $\ell$ is larger than half the supremum of all such radii. Smoothness at the pole ensures that such a radius exists. For $2\ell<t<1/(16A)$, no direction gives \eqref{eq:flatness} in $B_t$. Lemma~\ref{lem:nonflat} and $\Sigma\cap B_{At}\subset D$ give
\[
(1-s)\iint_{(\Sigma\cap B_{At})\times D}(1-\nu(x)\cdot\nu(y))K(x,y)\dd\sigma(y)\dd\sigma(x)\ge ct^{m-1-s}.
\]
For $2\ell<1/(16A)$, restrict the radius in \eqref{eq:radial-phi-basic} to $At$ and use Tonelli:
\[
\begin{aligned}
 &(1-s)\iint_{D^2}\Phi(x)(1-\nu(x)\cdot\nu(y))K\dd\sigma(y)\dd\sigma(x)\\
 &\ \ge c(1-s)\int_{2\ell}^{1/(16A)}(At)^{s-m}A
 \iint_{(D\cap B_{At})\times D}(1-\nu(x)\cdot\nu(y))K\dd\sigma(y)\dd\sigma(x)\dd t\\
 &\ \ge c(m)\int_{2\ell}^{1/(16A)} t^{s-m}t^{m-1-s}\dd t .
\end{aligned}
\]
If the interval is empty the integrand on $D^2$ is still nonnegative. Hence
\begin{equation}\label{eq:scale-payment}
(1-s)\iint_{D\times D}\Phi(x)\,(1-\nu(x)\cdot\nu(y))K(x,y)\dd\sigma(y)\dd\sigma(x)\ge c(m)\max\!\left\{0,\log\frac{1}{32A\ell}\right\}.
\end{equation}

By area growth, $0\le g\le e^2-1$, and \eqref{eq:scale-payment},
\[
\begin{aligned}
 \int_{D\cap B_{c_1\ell}}g(x)|x|^{-m}(|x|/\ell)^{\gamma/(m+1)}\dd\sigma(x)
 &\le C\sum_{j\ge0}2^{-j\gamma/(m+1)}\le C,\\
 \int_{D\cap(B_{1/4}\setminus B_{c_1\ell})}g(x)|x|^{-m}\dd\sigma(x)
 &\le C(1+|\log\ell|)\\
 &\le C+C(1-s)\iint_{D^2}\Phi(x)(1-\nu(x)\cdot\nu(y))K\dd\sigma(y)\dd\sigma(x).
\end{aligned}
\]
Multiply \eqref{eq:bad-defect}--\eqref{eq:outer-defect} by $g(x)$ and integrate. The two preceding bounds occur with the factor $(1-s)q^{-s}$. Thus
\begin{equation}\label{eq:integrated-defect}
\begin{aligned}
&-(1-s)\int_D g(x)\PV\int_D \frac{1+\nu(x)\cdot\nu(y)}{2}\bigl(\Phi(x)-\Phi(y)\bigr)K(x,y)\dd\sigma(y)\dd\sigma(x)\\
&\quad\le C\int_D g\dd\sigma+\frac q4 (1-s)\iint_{D\times D}g(x)\Phi(x)\,(1-\nu(x)\cdot\nu(y))K(x,y)\dd\sigma(y)\dd\sigma(x)\\
&\qquad +C(1-s)q^{-s}\left[1+(1-s)\iint_{D\times D}\Phi(x)\,(1-\nu(x)\cdot\nu(y))K(x,y)\dd\sigma(y)\dd\sigma(x)\right].
\end{aligned}
\end{equation}

Multiply \eqref{eq:sym-jacobi} by $\Phi(x)\ge0$ and integrate over $D$. Since $\Phi=0$ for $|x|\ge1$, only $x\in B_1$ enters \eqref{eq:sym-jacobi}. Formula \eqref{eq:radial-dirac}, with $f=g$, gives $g(0)$ when we move the integral operator onto $\Phi$. The right-hand side is bounded by $Cq(1-s)$ since $1+g\le e^2$ and $\int_D\Phi\,\dd\sigma\le C$. Thus
\[
\begin{aligned}
&g(0)+(1-s)\int_D g(x)\PV\int_D \frac{1+\nu(x)\cdot\nu(y)}{2}\bigl(\Phi(x)-\Phi(y)\bigr)K(x,y)\dd\sigma(y)\dd\sigma(x)\\
&\quad+\frac q2 (1-s)\iint_{D\times D}(1+g(x))\Phi(x)\,(1-\nu(x)\cdot\nu(y))K(x,y)\dd\sigma(y)\dd\sigma(x)\\
&\quad\le Cq(1-s).
\end{aligned}
\]
Substitution of \eqref{eq:integrated-defect} gives
\[
\begin{aligned}
 g(0)&+(1-s)\iint_{D^2}\left[\frac q2+\frac q4 g(x)-C(1-s)q^{-s}\right]
   \Phi(x)(1-\nu(x)\cdot\nu(y))K\dd\sigma(y)\dd\sigma(x)\\
 &\le C\int_D g\dd\sigma+C(1-s)(q+q^{-s}).
\end{aligned}
\]
By $(1-s)q^{-1-s}\le c(m)$ the bracket is at least $q/4$. Discard its nonnegative integral and rescale to obtain \eqref{eq:radial-mean}.
\end{proof}

\subsection{The contact bound and the three tests}

To verify the smallness condition on $(1-s)q^{-1-s}$ in Lemma~\ref{lem:radial-mean}, we first bound $\omega$ from below. We apply \eqref{eq:jacobi} at a minimum of $\omega/\eta$ and use Lemma~\ref{lem:moment} to control the linear term.

\begin{lemma}\label{lem:coarse}
Under the hypotheses of Theorem~\ref{thm:gradient}, for every base point $x'$ with $B'_d(x')\subset B'_{2R}$, the vertical normal at $x=(x',u(x'))$ satisfies
\begin{equation}\label{eq:coarse}
\omega(x)\ge c(m)(1-s)(1+\osc_{B'_{2R}}u/d)^{-m-2}.
\end{equation}
\end{lemma}

\begin{proof}
Scale lengths by $d/4$ and translate the base center to zero. On $\overline{B'_1}$ smoothness gives $\omega\ge(1+\max_{\overline{B'_1}}|Du|^2)^{-1/2}>0$. Hence $\omega/\eta$ tends to infinity as $\eta\downarrow0$ and attains its minimum at a graph point $x_*$ with $\eta(x_*')>0$. Only in this proof, write $\eta_*=\eta(x_*')$, $\omega_*=\omega(x_*)$ and $K=K(x_*,y)$. The minimum gives $\omega(y)\ge(\omega_*/\eta_*)\eta(y')$. Choose $\rho=c\sqrt{\eta_*}$ with $c$ small. Taylor expansion on the base, \eqref{eq:moment} and \eqref{eq:moments-upper} give
\[
\begin{aligned}
 &(1-s)\PV\int_{\Sigma\cap C_2}\eta((y-x_*)/\rho)
     [\eta_*\nu(x_*)\cdot\nu(y)-\eta(y')]K\dd\sigma(y)\\
 &\le C(1-s)|D\eta(x_*')|\rho^{-s}
   +C(1-s)\int_{\Sigma\cap C_2\cap B_\rho(x_*)}|y-x_*|^2 K\dd\sigma(y)\\
 &\quad+[C\rho|D\eta(x_*')|-\eta_*](1-s)
 \int_{\Sigma\cap C_2}\eta((y-x_*)/\rho)(1-\nu(x_*)\cdot\nu(y))K\dd\sigma(y)\\
 &\le C(1-s)\sqrt{\eta_*}\rho^{-s}+C\rho^{1-s}\le C .
\end{aligned}
\]
Indeed $|D\eta(x_*')|\le C\sqrt{\eta_*}$, so the bracket is $\le-\eta_*/2$, and $\sqrt{\eta_*}\rho^{-s}=c^{-1}\rho^{1-s}$. Splitting \eqref{eq:jacobi} at $x_*$ now gives
\[
\begin{aligned}
 &(1-s)\int_{\Sigma\cap C_2}(1-\eta((y-x_*)/\rho))\omega(y)K\dd\sigma(y)\\
 &\le\frac{\omega_*}{\eta_*}(1-s)\PV\int_{\Sigma\cap C_2}
 \eta((y-x_*)/\rho)[\eta_*\nu(x_*)\cdot\nu(y)-\eta(y')]K\dd\sigma(y)\\
 &\quad+(1-s)\omega_*\int_{\Sigma\cap C_2}
 (1-\eta((y-x_*)/\rho))\nu(x_*)\cdot\nu(y)K\dd\sigma(y)+C(1-s)\omega_*\\
 &\le C\frac{\omega_*}{\eta_*}
   \bigl[1+(1-s)\eta_*\rho^{-1-s}\bigr]\le C\frac{\omega_*}{\eta_*}.
\end{aligned}
\]
The far integral is at most $\int_{(\Sigma\cap C_2)\setminus B_{3\rho/4}(x_*)}K\dd\sigma\le C\rho^{-1-s}$ by \eqref{eq:moments-upper}, and $\eta_*\rho^{-1-s}=c^{-2}\rho^{1-s}$. Over $3/2<|y'|<2$ the far cutoff is one. Since $\omega(y)\dd\sigma=\dd y'$ and $\eta(0)=1$,
\[
 \omega(0,u(0))\ge\frac{\omega_*}{\eta_*}
 \ge c(1-s)\int_{B'_2\setminus B'_{3/2}}\frac{\dd y'}{|(y',u(y'))-x_*|^{m+1+s}}
 \ge c(1-s)(1+\osc_{B'_2}u)^{-m-2}.
\]
In the original coordinates the rescaled oscillation is at most $4\osc_{B'_{2R}}u/d$, giving \eqref{eq:coarse}.
\end{proof}

To bound the surface integral in \eqref{eq:radial-mean}, we use the graph equation, the Jacobi inequality and the graph equation in succession, following the order of the three tests in Trudinger \cite[(23)--(27)]{Trudinger}.

\begin{lemma}\label{lem:height}
Let $P$ lie on the graph of $u$ and let $P'$ be its projection. Take $R>40r$ and assume that $C_{2R}(P')$ lies in the smooth minimizing region, so \eqref{eq:jacobi} holds above $B'_{20r}(P')$. If $0<q\le1/2$, $\omega^{-q}\le e^2$ on $\Sigma\cap C_R(P')$, and $(1-s)e^{1/(2q)}\le1$, then
\begin{equation}\label{eq:height}
\int_{\Sigma\cap B_{4r}(P)}(\omega^{-q}-1)\dd\sigma\le C(m)q r^{m-1}\bigl(r+\osc_{B'_{20r}(P')}u\bigr).
\end{equation}
\end{lemma}

\begin{proof}
Translate $P$ to zero and divide all lengths by $r$, including $R$, so $u(0)=0$. The two base cutoffs $\eta(x'/8)$ and $\eta(x'/16)$ stand for $\eta((x'-P')/(8r))$ and $\eta((x'-P')/(16r))$ in the original coordinates. The first equals one on $B'_6$. The second equals one on $B'_{12}$ and thus on the entire support of the first.

Set $g=\omega^{-q}-1$. The assumptions give $\omega^{-1}\le e^{2/q}\le(1-s)^{-4}$ on $\Sigma\cap C_R$, and $1-s\le e^{-1/(2q)}\le q$ since $e^t\ge 2t$ for $t>0$. Since $\omega^{-1}\le(1-s)^{-4}$, the graph area formula gives $|B'_{12}|\le\int_{\Sigma\cap C_R}\eta(x'/16)\dd\sigma(x)\le C(1-s)^{-4}$.

The weak graph equation of Cozzi and Lombardini \cite[(1.5)--(1.6), Theorem~1.10(iv)$\Rightarrow$(ii)]{CL}, multiplied by $1-s$, is
\begin{equation}\label{eq:weak-graph}
\begin{gathered}
(1-s)\iint_{\R^m\times\R^m}
F_s\!\left(\frac{u(x')-u(y')}{|x'-y'|}\right)
\frac{\varphi(x')-\varphi(y')}{|x'-y'|^{m+s}}\dd x'\dd y'=0,\\
F_s(t)=\int_0^t(1+\tau^2)^{-(m+1+s)/2}\dd\tau .
\end{gathered}
\end{equation}
The test $\varphi$ is compactly supported in $B'_{20}$. Lipschitz tests follow by smooth approximation, since $F_s$ is uniformly bounded and the kernel acting on a Lipschitz difference is integrable at fixed $s<1$.
For $0<h<(1-s)/100$ take the surface average
\[
 g_h(x)=\frac{\int_{\Sigma\cap C_R}\eta((y-x)/h)\,g(y)\,\dd\sigma(y)}
                 {\int_{\Sigma\cap C_R}\eta((y-x)/h)\,\dd\sigma(y)} .
\]
The center $x$ may vary in an ambient neighborhood of the graph above $B'_8$. For the first test in \eqref{eq:weak-graph}, take
\[
 \varphi(x')=\eta(x'/8)\,g_h(x',u(x'))\,\min\{(u(x')+8)_+,16\}.
\]
The Jacobi test below selects a radius $h$ at which the derivative of this weight is controlled. Both graph tests have the form $\varphi=f\,\psi(u)$ with $f\ge0$ and $\psi$ bounded and nondecreasing. The test difference equals
\[
f(x')\bigl(\psi(u(x'))-\psi(u(y'))\bigr)+\psi(u(y'))\bigl(f(x')-f(y')\bigr).
\]
Oddness and monotonicity of $F_s$ make the first term, including its factor $F_s$, nonnegative. Using the equation to estimate it by the second term gives
\begin{equation}\label{eq:height-diffs}
\begin{split}
&(1-s)\iint_{\R^m\times\R^m}
f(x') F_s\!\left(\frac{u(x')-u(y')}{|x'-y'|}\right)
\frac{\psi(u(x'))-\psi(u(y'))}{|x'-y'|^{m+s}}\dd x'\dd y'
\\
&\hspace{25mm}\le C \|\psi\|_\infty\bigl(\|f\|_{L^1(\R^m)}+\|Df\|_{L^1(\R^m)}\bigr).
\end{split}
\end{equation}
Indeed the integral in $x'$ of $|f(x'+z)-f(x')|$ is bounded by $|z|\|Df\|_1$ for $|z|\le1$ and by $2\|f\|_1$ for $|z|>1$. Their radial integrals are $(1-s)\int_0^1\rho^{-s}\dd\rho=1$ and $(1-s)\int_1^\infty\rho^{-1-s}\dd\rho=(1-s)/s$. To get a weighted area bound from \eqref{eq:height-diffs}, we need both a derivative bound for $\eta(x'/8) g_h$ and a comparison between the positive height differences and $\omega|Du|^2$ on the base.

The second test in \eqref{eq:jacobi} is $\eta(x'/16)^2/\omega(x)$, as in Cabr\'e--Cozzi \cite[proof of Lemma~6.6, (6.29)]{CC}. By symmetry,
\[
\begin{aligned}
 &\int_{\Sigma\cap C_R}\frac{\eta(x'/16)^2}{\omega(x)}
      L_{\Sigma\cap C_R}\omega(x)\dd\sigma(x)\\
 &=\frac{1-s}{2}\iint_{(\Sigma\cap C_R)^2}(\omega(x)-\omega(y))
 \left(\frac{\eta(x'/16)^2}{\omega(x)}-\frac{\eta(y'/16)^2}{\omega(y)}\right)
 K(x,y)\dd\sigma(y)\dd\sigma(x).
\end{aligned}
\]
Expanding a square gives
\[
\begin{aligned}
 &(\omega(x)-\omega(y))
 \left(\frac{\eta(x'/16)^2}{\omega(x)}-\frac{\eta(y'/16)^2}{\omega(y)}\right)\\
 &\quad=(\eta(x'/16)-\eta(y'/16))^2
 -\frac{(\eta(x'/16)\omega(y)-\eta(y'/16)\omega(x))^2}{\omega(x)\omega(y)} .
\end{aligned}
\]
Restrict the square to the graph over $B'_{12}$, where $\eta(x'/16)=1$. This term and the geometric term in \eqref{eq:jacobi} therefore have upper bound $C\int_{\Sigma\cap C_R}\eta(x'/16)\dd\sigma$. To see the bound, use
\[
|\eta(x'/16)-\eta(y'/16)|^2\le C(\eta(x'/16)+\eta(y'/16))\min\{1,|x-y|^2\},
\]
which follows since $\eta(x'/16)$ is a square. Symmetry and \eqref{eq:moments-upper} bound the contribution of this cutoff difference by $C\int_{\Sigma\cap C_R}\eta(x'/16)\dd\sigma$. The truncation error is bounded by $C(1-s)\int_{\Sigma\cap C_R}\eta(x'/16)^2\dd\sigma\le C\int_{\Sigma\cap C_R}\eta(x'/16)\dd\sigma$. Consequently
\begin{equation}\label{eq:height-normal}
\begin{split}
&(1-s)\iint_{(\Sigma\cap C_{12})^2}
 \frac{|\omega(x)-\omega(y)|^2}{\omega(x)\omega(y)}
 K(x,y)\dd\sigma(y)\dd\sigma(x)\\
&\quad+(1-s)\iint_{(\Sigma\cap C_R)^2}\eta(x'/16)^2|\nu(x)-\nu(y)|^2 K(x,y)\dd\sigma(y)\dd\sigma(x)\le C\int_{\Sigma\cap C_R}\eta(x'/16)\dd\sigma.
\end{split}
\end{equation}
By symmetry take $\omega(x)\ge\omega(y)$ and $t=\omega(x)/\omega(y)$. Since $1-t^{-q}\le q\log t\le q(t-1)/\sqrt t$ and $\omega(y)^{-q}\le e^2$,
\[
\begin{aligned}
 |g(x)-g(y)|&=\omega(y)^{-q}(1-t^{-q})
 \le e^2 q\,\frac{|\omega(x)-\omega(y)|}{\sqrt{\omega(x)\omega(y)}},\\
 (1-s)\iint_{(\Sigma\cap C_{12})^2}|g(x)-g(y)|^2K\dd\sigma(y)\dd\sigma(x)
 &\le Cq^2\int_{\Sigma\cap C_R}\eta(x'/16)\dd\sigma(x)
\end{aligned}
\]
by \eqref{eq:height-normal}.

We now estimate the surface average in the first multiplier. For $x\in\Sigma\cap C_8$, Lemma~\ref{lem:geometry} bounds the area of $\Sigma\cap B_h(x)$ above by $Ch^m$ and that of $\Sigma\cap B_{3h/4}(x)$ below by $ch^m$. Thus the denominator in $g_h$ is comparable to $h^m$. Differentiate its quotient in the ambient center $x$. Subtracting the constant $g(x)$ from the integrand and using $|D_x[\eta((y-x)/h)]|\le C/h$ gives
\[
|g_h(x)-g(x)|+h|Dg_h(x)|\le C h^{-m}
\int_{\Sigma\cap B_h(x)}|g(y)-g(x)|\dd\sigma(y).
\]
Square by Cauchy--Schwarz and integrate on $\Sigma\cap C_{8}$. For the first term use $h^{-m}\le h^{1+s}K(x,y)$ on $|x-y|<h$. For the gradient term, integrate the bound in $h$ against $(1-s)h^{-s}\dd h$ and use $\int_{|x-y|}^{(1-s)/100} h^{-m-2-s}\dd h\le C K(x,y)$. This gives
\[
\sup_{h<(1-s)/100}\|g_h-g\|_{L^2(\Sigma\cap C_{8})}^2
+(1-s)\int_0^{(1-s)/100} h^{-s}\|Dg_h\|_{L^2(\Sigma\cap C_{8})}^2\dd h\le Cq^2\int_{\Sigma\cap C_R}\eta(x'/16)\dd\sigma.
\]
The averaging weight has mass $((1-s)/100)^{1-s}\ge c>0$. Select $h$ for which the two squared norms are bounded by $Cq^2\int_{\Sigma\cap C_R}\eta(x'/16)\dd\sigma$. On the base use $|D_{x'}[g_h(x',u(x'))]|\le |Dg_h|/\omega$, $\dd x'=\omega\dd\sigma$ and $0\le\omega g=\omega^{1-q}-\omega\le q$ to get
\[
\begin{aligned}
 \int_{B'_8}g_h\dd x'&\le q\HH^m(\Sigma\cap C_8)
          +\int_{\Sigma\cap C_8}|g_h-g|\dd\sigma,\\
 \int_{B'_8}|D_{x'}[g_h(x',u(x'))]|\dd x'&\le\int_{\Sigma\cap C_8}|Dg_h|\dd\sigma.
\end{aligned}
\]
Cauchy--Schwarz, the product rule and $\eta(x'/16)=1$ on $C_8$ now give
\begin{equation}\label{eq:height-weight}
\begin{aligned}
 &\|\eta(x'/8)g_h\|_{L^1(\R^m)}
 +\|D(\eta(x'/8)g_h)\|_{L^1(\R^m)}
 +\int_{\Sigma\cap C_{8}}|g_h-g|\dd\sigma \\
 &\qquad\le Cq\int_{\Sigma\cap C_R}\eta(x'/16)\dd\sigma .
\end{aligned}
\end{equation}
In the first two norms $g_h$ is evaluated on the graph.

To compare the height difference quotients with $Du(x')\cdot\theta$, we first bound $Du(x')-Du(y')$ using the integral of the squared normal difference in \eqref{eq:height-normal}. On the graph $Du=-\nu'/\omega$, hence
$|Du(x')-Du(y')|\le C(1-s)^{-8}|\nu(x)-\nu(y)|$. The chord satisfies $|x-y|\le(1-s)^{-4}|x'-y'|$ for endpoints above $B'_{20}$, by integrating $\omega^{-1}$ along the segment on the base. Since $\dd x'=\omega\dd\sigma\le\dd\sigma$, \eqref{eq:height-normal} bounds the following integral from above by $C(1-s)^{-4(m+5+s)}\int_{\Sigma\cap C_R}\eta(x'/16)\dd\sigma$. Using $\int_{\Sigma\cap C_R}\eta(x'/16)\dd\sigma\le C(1-s)^{-4}$ and $s<1$ gives
\begin{equation}\label{eq:height-grad-diff}
(1-s)\iint_{B'_{16}\times B'_{20}}
\frac{\eta(x'/16)^2 |Du(x')-Du(y')|^2}{|x'-y'|^{m+1+s}}\dd x'\dd y'
\le C(1-s)^{-4(m+7)}.
\end{equation}
Take $\delta=(1-s)^{2m+15}$. With $x'\in B'_{16}$, $0<\rho<\delta$, and $|\theta|=1$, the fundamental theorem of calculus gives
\[
 \left|\frac{u(x'+\rho\theta)-u(x')}{\rho}-Du(x')\cdot\theta\right|^2
 \le\int_0^1|Du(x'+\tau\rho\theta)-Du(x')|^2\dd\tau .
\]
Multiply by $(1-s)\eta(x'/16)^2\rho^{-s}$ and integrate in $x',\theta,\rho$. Changing variables to $v=\tau\rho$ yields the upper bound
\[
\begin{aligned}
 &(1-s)\int_0^1\tau^{s-1}\int_{B'_{16}}\int_{\mathbb S^{m-1}}\int_0^{\tau\delta}
 \eta(x'/16)^2|Du(x'+v\theta)-Du(x')|^2v^{-s}\dd v\dd\theta\dd x'\dd\tau\\
 &\le\frac{\delta^2}{s}(1-s)\iint_{B'_{16}\times B'_{20}}
 \eta(x'/16)^2\frac{|Du(y')-Du(x')|^2}{|y'-x'|^{m+1+s}}\dd y'\dd x'.
\end{aligned}
\]
Use \eqref{eq:height-grad-diff} and Cauchy--Schwarz to obtain
\begin{equation}\label{eq:height-diff-error}
\begin{aligned}
&(1-s)\int_{B'_{16}}\int_{\mathbb S^{m-1}}\int_0^\delta
\eta(x'/16)\,\left|\frac{u(x'+\rho\theta)-u(x')}{\rho}-Du(x')\cdot\theta\right|
\frac{\dd\rho}{\rho^s}\dd\theta\dd x'\\
&\qquad\le C\delta(1-s)^{-2(m+7)}=C(1-s).
\end{aligned}
\end{equation}
The radial integral over $0<\rho<\delta$ has a positive lower bound depending only on $m$. Indeed,
\[
 (1-s)\int_0^\delta \rho^{-s}\,\dd\rho
 =\delta^{1-s}=\exp\bigl[(2m+15)(1-s)\log(1-s)\bigr]\ge e^{-(2m+15)/e}.
\]
Thus \eqref{eq:height-diff-error} makes the replacement error small while leaving a uniform positive coefficient in the area estimate.

For every $p\in\R^m$, $tF_s(t)\ge tF_1(t)$ and
\[
\int_{\mathbb S^{m-1}}(p\cdot\theta)F_1(p\cdot\theta)\dd\theta
=|p|^2\int_{B'_1}\frac{\dd z}{(1+(p\cdot z)^2)^{(m+2)/2}}
=|B'_1|\frac{|p|^2}{\sqrt{1+|p|^2}} .
\]
The first equality is the divergence theorem on the ball. For the second, integrate slices perpendicular to $p$. The substitution $t\mapsto t\sqrt{(1+|p|^2)/(1+|p|^2t^2)}$ gives the stated value.

Let $U\subset B'_{16}$ be a base region on which the height truncation preserves every short difference: for $x'\in U$, $0<\rho<\delta$ and $|\theta|=1$,
\[
 \psi(u(x'+\rho\theta))-\psi(u(x'))=u(x'+\rho\theta)-u(x').
\]
Restrict the nonnegative integral in \eqref{eq:height-diffs} to these pairs. For $0\le f\le(e^2-1)\eta(x'/16)$ the function $t\mapsto tF_s(t)$ has uniformly bounded derivative. By \eqref{eq:height-diff-error}, the angular identity above and $\delta^{1-s}\ge c(m)$,
\begin{equation}\label{eq:height-diff-area}
\begin{aligned}
 &(1-s)\int_U\int_{|z|<\delta}f(x')
 \frac{u(x'+z)-u(x')}{|z|}F_s\!\left(\frac{u(x'+z)-u(x')}{|z|}\right)
 \frac{\dd z\dd x'}{|z|^{m+s-1}}\\
 &\quad\ge \delta^{1-s}\int_U f(x')\int_{\mathbb S^{m-1}}
 (Du(x')\cdot\theta)F_1(Du(x')\cdot\theta)\dd\theta\dd x'-C(1-s)\\
 &\quad\ge c(m)\int_U f\,\omega |Du|^2\dd x'-C(1-s).
\end{aligned}
\end{equation}
Return to the first test, with $U=\{x':|x'|^2+u(x')^2<16\}$. For $x'\in U$, $|u(x')|<4$ and $|u(x'+z)-u(x')|\le(1-s)^{-4}\delta<1/2$, so both heights are in the affine range of $\psi$. On $U$ one has $f=g_h$. Equations \eqref{eq:height-diffs}, \eqref{eq:height-weight} and \eqref{eq:height-diff-area} give
\[
 \int_U g_h\omega|Du|^2\dd x'\le C(1-s)+Cq\int_{\Sigma\cap C_R}\eta(x'/16)\dd\sigma
 \le Cq\int_{\Sigma\cap C_R}\eta(x'/16)\dd\sigma .
\]
Here $1-s\le q$. Using $\omega^{-1}\le1+\omega|Du|^2$ and \eqref{eq:height-weight},
\begin{equation}\label{eq:height-weighted-area}
\begin{aligned}
 \int_{\Sigma\cap B_4}g\dd\sigma
 &\le\int_{\Sigma\cap C_8}|g-g_h|\dd\sigma
   +\int_U g_h\dd x'+\int_U g_h\omega|Du|^2\dd x'\\
 &\le Cq\int_{\Sigma\cap C_R}\eta(x'/16)\dd\sigma .
\end{aligned}
\end{equation}

Finally, take $\varphi(x')=\eta(x'/16)u(x')$ in \eqref{eq:weak-graph}. Since $u(0)=0$, one has $|u|\le\osc_{B'_{20}}u$ on $B'_{20}$. To apply \eqref{eq:height-diffs}, take $\psi$ equal to the identity on $[-\osc_{B'_{20}}u,\osc_{B'_{20}}u]$ and extend it by constants on both sides. Then $\eta(x'/16)\psi(u)=\varphi$ everywhere, including outside $B'_{16}$, and $\|\psi\|_\infty=\osc_{B'_{20}}u$. For $x'\in B'_{16}$ and $|z|<\delta$, both $x'$ and $x'+z$ belong to $B'_{20}$, so $\psi$ preserves their height difference. Taking $U=B'_{16}$ and $f(x')=\eta(x'/16)$, \eqref{eq:height-diff-area} gives the lower bound below. By nonnegativity, \eqref{eq:height-diffs} gives the upper bound, since $\|f\|_1+\|Df\|_1\le C(m)$. Thus
\[
\begin{aligned}
&c(m)\int_{B'_{16}}\eta(x'/16)\,\omega |Du|^2\,\dd x'-C(1-s)\\
&\quad\le (1-s)\int_{B'_{16}}\int_{|z|<\delta}
  \eta(x'/16)\frac{u(x'+z)-u(x')}{|z|}
  F_s\!\left(\frac{u(x'+z)-u(x')}{|z|}\right)
  \frac{\dd z\,\dd x'}{|z|^{m+s-1}}\\
&\quad\le C(m)\osc_{B'_{20}}u .
\end{aligned}
\]
Using $\omega^{-1}-1\le\omega|Du|^2$ gives
\[
 \int_{\Sigma\cap C_R}\eta(x'/16)\,\dd\sigma
 \le |B'_{16}|+\int_{B'_{16}}\eta(x'/16)\omega|Du|^2\,\dd x'
 \le C(m)(1+\osc_{B'_{20}}u).
\]
Returning to the original coordinates in \eqref{eq:height-weighted-area} gives
\[
 \int_{\Sigma\cap B_{4r}(P)}g\dd\sigma
 \le C q r^m\left(1+\frac{\osc_{B'_{20r}(P')}u}{r}\right).
\]
\end{proof}

We now choose $P,q,r$ so that $g(P)=e-1$ and both local estimates apply.

\begin{proof}[Proof of Theorem~\ref{thm:gradient}]
On the closed ball $\overline{B'_{3R/2}}$, take
\[
 \max_{|x'|\le 3R/2}
       \bigl(\sqrt{1+|Du(x')|^2}\bigr)^{1-2|x'|/(3R)} .
\]
The exponent equals $\operatorname{dist}(x',\partial B'_{3R/2})/(3R/2)$. Smoothness inside $B'_{2R}$ makes the displayed function continuous on the closed auxiliary ball. It equals one on its boundary. On the target ball $B'_R$ the exponent is at least $1/3$, so the cube of the maximum bounds $\sqrt{1+|Du|^2}$ there.

Fix a dimensional $C_0\ge2$, to be chosen at the end, and suppose the maximum exceeds $\exp[C_0(1+\osc_{B'_{2R}}u/R)]$. It is then attained at an interior point $P'$. Let $P=(P',u(P'))$, $d=3R/2-|P'|$ and choose $q>0$ by $\omega(P)^{-q}=e$. At the maximum we have
\[
 \exp\!\left(\frac{2d}{3Rq}\right)>\exp[C_0(1+\osc_{B'_{2R}}u/R)],
 \qquad \frac1q>\frac{3C_0(R+\osc_{B'_{2R}}u)}{2d}\ge C_0(1+\osc_{B'_{2R}}u/d),
\]
where $d\le3R/2$. In particular $q\le1/2$.

Maximality also bounds the weight in a neighborhood of $P$. If $|x'-P'|<d/2$, the triangle inequality gives $3R/2-|x'|\ge d/2$. For the corresponding graph point $x=(x',u(x'))$, we have $\omega(x)^{-1}\ge1$, so
\[
 (\omega(x)^{-1})^{d/(3R)}
 \le (\omega(x)^{-1})^{1-2|x'|/(3R)}
 \le (\omega(P)^{-1})^{2d/(3R)} .
\]
Raising to the power $3Rq/d$ gives $\omega(x)^{-q}\le\omega(P)^{-2q}=e^2$. Thus
\begin{equation}\label{eq:radial-contact}
g=\omega^{-q}-1,\qquad g(P)=e-1,\qquad 0\le g\le e^2-1\quad\text{on }\Sigma\cap C_{d/2}(P').
\end{equation}

Since $B'_d(P')\subset B'_{3R/2}\Subset B'_{2R}$, the preliminary estimate \eqref{eq:coarse} gives $e^{1/q}\le C(1-s)^{-1}(1+\osc_{B'_{2R}}u/d)^{m+2}$. Choose the dimensional $C_0(m)$ large enough that $Ct^{m+2}\le e^{C_0t/2}$ for every $t\ge1$. At the maximum this gives
\[
 (1-s)e^{1/q}\le C(1+\osc_{B'_{2R}}u/d)^{m+2}
 \le e^{C_0(1+\osc_{B'_{2R}}u/d)/2}\le e^{1/(2q)} .
\]
Hence $(1-s)e^{1/(2q)}\le1$. Since $q\le1/2$, we have $q^{-1-s}\le q^{-2}$. In particular,
\begin{equation}\label{eq:q-regime}
(1-s)q^{-1-s}\le (1-s)^{1/2}q^{-2}e^{-1/(4q)}\le C(m)(1-s)^{1/2}.
\end{equation}

Take $r=d/200$ and $D=\Sigma\cap B_{4r}(P)$. Lemma~\ref{lem:radial-mean} applies in $C_{40r}(P')$. In Lemma~\ref{lem:height} choose the outer radius $d/4$. Both choices are admissible since
\[
 20r=d/10<40r=d/5<d/4,\qquad
 C_{40r}(P')\subset C_{d/4}(P')\subset C_{d/2}(P')\subset C_{2R}.
\]
By \eqref{eq:radial-contact}, $\omega^{-q}\le e^2$ throughout $\Sigma\cap C_{d/2}(P')$.

By \eqref{eq:q-regime} we may choose $s_*(m)$ so that $(1-s)q^{-1-s}\le c(m)$, the hypothesis of Lemma~\ref{lem:radial-mean}. Insert \eqref{eq:height} into \eqref{eq:radial-mean} on this same $D$, use $g(P)=e-1$, and divide by $q$. Since $r=d/200$ and $\osc_{B'_{20r}(P')}u\le \osc_{B'_{2R}}u$, we obtain
\[
 \frac{e-1}{q}\le C(m)\left(1+\frac{\osc_{B'_{2R}}u}{r}\right)+C(m)(1-s)q^{-1-s}
                 \le C(m)\left(1+\frac{\osc_{B'_{2R}}u}{d}\right).
\]
The last inequality uses \eqref{eq:q-regime}. All constants on the right are independent of $C_0$. We can therefore choose $C_0=C_0(m)$ large enough to contradict $q^{-1}>C_0(1+\osc_{B'_{2R}}u/d)$. The weighted maximum is at most $\exp[C_0(1+\osc_{B'_{2R}}u/R)]$. Its cube bounds $\omega^{-1}$ on $B'_R$, proving \eqref{eq:main}.
\end{proof}

\section*{Declaration on the use of generative AI}

The author supplied the main ideas and directed the proof strategy. The author is responsible for the mathematical arguments and the manuscript. ChatGPT (OpenAI) was used to assist in exploring individual proof steps, checking algebra and scaling, editing the LaTeX source, and preparing the translation.

%--------------------------------------------------------
% Bibliography
% ------------------------------------------------
\small
\bibliographystyle{alpha}
\bibliography{refs}
\medskip
\noindent
(Shuning Xu) Department of Mathematics, University of Science and Technology of China, Hefei, 230026, Anhui Province, China.\;
Email address: xushuning@mail.ustc.edu.cn

\end{document}